\documentclass{amsart}

\usepackage{a4wide} 
\usepackage{amsmath,amssymb,amsthm}
\usepackage{enumerate} 
\usepackage{bbm}

\newtheorem{theorem}{Theorem}[section]
\newtheorem{lemma}[theorem]{Lemma}
\newtheorem{proposition}[theorem]{Proposition}
 
\theoremstyle{definition}
\newtheorem{definition}[theorem]{Definition}

\theoremstyle{remark}
\newtheorem{remark}[theorem]{Remark}

\numberwithin{equation}{section}

\newcommand{\Ann}[1]{{#1}^0} 
\newcommand{\e}{{\rm e}}
\newcommand{\ep}{{\rm e}} 
\newcommand{\f}{{\rm f}}
\newcommand{\ccdot}{\!\cdot\!}  
\newcommand{\hook}{\lrcorner \,}
\newcommand{\id}{\mathbbm{1}}

\newcommand{\bR}{\mathbb{R}}

\renewcommand{\i}{{\rm i}}
\newcommand{\mm}[4]{\left(\begin{matrix}#1&#2\\#3&#4\end{matrix}\right)}
\newcommand{\mmm}[9]{\left(\begin{matrix}#1&#2&#3\\#4&#5&#6\\#7&#8&#9\end{matrix}\right)}

\renewcommand{\Re}{{\rm Re}} 
\renewcommand{\Im}{{\rm Im}}
 
\newcommand{\im}{{\rm im}}
\renewcommand{\dim}{{\rm dim\ }}

\newcommand{\solv}{\mathfrak {r}}
\renewcommand{\sl}{\mathfrak {sl}} 
\newcommand{\su}{\mathfrak {su}}
\newcommand{\g}{\mathfrak{g}} 
\newcommand{\h}{\mathfrak{h}}
\newcommand{\efr}{\mathfrak{e}}

\newcommand{\GL}{{\rm GL}} 
\newcommand{\SL}{{\rm SL}}
\newcommand{\U}{{\rm U}} 
\newcommand{\SU}{{\rm SU}}
\newcommand{\G}{{\rm G}}

\newcommand{\op}{\oplus} 
\newcommand{\ot}{\otimes}
\renewcommand{\o}{\omega} 
\renewcommand{\a}{\alpha}
\renewcommand{\b}{\beta} 
\newcommand{\la}{\lambda}
 
\renewcommand{\L}{\Lambda}
\newcommand{\s}{\sigma} 
\renewcommand{\^}{\wedge}

\begin{document}

\title[Half-flat structures on products of three-dimensional Lie groups]{Half-flat structures on products \\of three-dimensional Lie groups}

\author{Fabian Schulte-Hengesbach}

\address{Department Mathematik, Universit\"at Hamburg,
  Bundesstra{\ss}e 55, D-20146 Hamburg, Germany}
\email{schulte-hengesbach@math.uni-hamburg.de}

\thanks{The results of this article will be part of the author's
  PhD-thesis supported by the SFB 676 ``Particles, strings and the
  early universe: the geometry of matter and space-time'' of the
  Deutsche Forschungsgemeinschaft.}



\begin{abstract}
  We classify six-dimensional Lie groups which admit a left-invariant
  half-flat $\SU(3)$-structure and which split in a direct product of
  three-dimensional factors. Moreover, a complete list of those direct
  products is obtained which admit a left-invariant half-flat
  $\SU(3)$-structure such that the three-dimensional factors are
  orthogonal. Similar classification results are proved for
  left-invariant half-flat $\SL(3,\bR)$-structures on direct products
  with either definite and orthogonal or isotropic factors.
\end{abstract}

\maketitle

\section{Introduction}
An $\SU(3)$-structure $(g,J,\o,\Psi)$ on a six-dimensional manifold
$M$ consists of a Riemannian metric $g$, an orthogonal almost complex
structure $J$, the fundamental two-form $\omega = g(.\,,J\,.)$ and a
complex-valued $(3,0)$-form $\Psi$ of constant length. If furthermore
the exterior system
\[ d( \omega \wedge \omega) = 0 \;, \qquad d (\Re \Psi) = 0 \; ,\] is
satisfied, the $\SU(3)$-structure is called \emph{half-flat}. This
notion was introduced in \cite{ChSa} where $\SU(3)$-structures are
classified by irreducible components of the intrinsic torsion. The
main motivation for studying half-flat $\SU(3)$-structures is their
close relation to parallel $\G_2$-structure via the Hitchin flow. On
the one hand, a parallel $\G_2$-structure on a seven-manifold induces
a half-flat $\SU(3)$-structure on every oriented hypersurface. On the
other hand, half-flat $\SU(3)$-structures on a compact six-manifold
$M$ can be embedded in a seven-manifold with parallel $\G_2$-structure
as follows. Given a (global) solution of the Hitchin flow on an
interval $I$ which is a half-flat $\SU(3)$-structure at a time
$t_0$, there is a parallel $\G_2$-structure on $M \times I$,
\cite{Hi1}. In fact, the proof is generalised to non-compact
six-manifolds in \cite{CLSS}.

Another motivation for the study of half-flat $\SU(3)$-structures
comes from string theory and supergravity which discusses them as
candidates for internal spaces of compactifications in the presence of
background fluxes (\cite{GLMW} or, more recently, \cite{GLM} and
references therein).

In the mathematical literature, half-flat $\SU(3)$-structures have
been studied intensively on nilmanifolds. For instance, a
classification under different additional assumptions is obtained in
\cite{CF}, \cite{ChSw} and \cite{CT}. Very recently, the
classification of nilmanifolds admitting invariant half-flat
$\SU(3)$-structures without any further restrictions has been obtained
in \cite{C}. Apart from the nilpotent case, examples and constructions
of half-flat $\SU(3)$-structures can be found in \cite{TV} and
\cite{AFFU}. The Ricci curvature of a half-flat $\SU(3)$-structure is
computed in \cite{BV} and \cite{AC}.

In this article, we ask the question which direct products of two
three-dimensional Lie groups admit a left-invariant half-flat
$\SU(3)$-structure. There are 12 isomorphism classes of
three-dimen\-sional Lie algebras (see tables \ref{unilist} and
\ref{nonunilist}), if we count the two Bianchi classes which depend on
a continuous parameter as three classes characterised by the property
that the parameter can be deformed continuously without leaving the
class. Thus, we have to consider $78=\binom{13}{2}$ classes of direct
sums in total after reducing the problem to the Lie algebra as usual.

Initially, we tried to find a classification by a direct proof which
avoids the verification of the existence or non-existence case by
case. However, this was only successful when we asked for the
existence of a half-flat $\SU(3)$-structure $(g,J,\o,\Psi)$ such that
the two factors are orthogonal with respect to the metric $g$. The
result is that exactly 15 classes admit such an $\SU(3)$-structure, 11
of which are unimodular and comply with a regular pattern, whereas the
remaining four do not seem to share many properties. Given, that the
additional assumption is rather strong and the proof, which is
presented in section \ref{orth}, is already quite technical, an answer
to the initial question with this method cannot be expected. However,
an advantage of the assumption of a Riemannian product is the fact
that the curvature is completely determined by the Ricci tensors of
the three-dimensional factors and that the possible Ricci tensors of
left-invariant metrics on three-dimensional Lie groups are classified
in \cite{Mi}. Furthermore, we remark that a basis is introduced in
Lemma \ref{u3ortho} which is well adapted to an almost Hermitian
structure on a Riemannian product of three-manifolds which could be
useful beyond the framework of half-flat structures.

A completely different method is used in \cite{C} for classifying the
nilmanifolds admitting an arbitrary half-flat $\SU(3)$-structure. An
obstruction to the existence of a half-flat $\SU(3)$-structure is
introduced in terms of the cohomology of a double complex which can be
constructed on most of the nilpotent Lie algebras. In our situation,
such a double complex can be constructed if and only if both Lie groups
are solvable. However, as the methods of homological algebra turn out
not to be advantageous for our problem, we prove a simplified version
of the obstruction condition in section \ref{sectionobst}. This
obstruction is applied directly to 41 isomorphism classes of direct
sums in section \ref{classi}. Two classes resist the obstruction,
although they do not admit a half-flat structure either, which is
proved individually by finding refined obstruction conditions. The
remaining 35 direct sums, including all unimodular direct sums and all
non-solvable direct sums, admit a half-flat $\SU(3)$-structure which is
proved by giving one explicit example in each case in the appendix. We
point out that the products of unimodular three-dimensional Lie groups
are particularly interesting since they admit co-compact lattices,
\cite{RV}.

In fact, the most time-consuming part of the classification was the
construction of examples of half-flat structures for the 20=35-15
classes which do not admit an ``orthogonal'' half-flat
$\SU(3)$-structure. The construction essentially relies on the fact
that a left-invariant half-flat $\SU(3)$-structure is defined by a
pair $(\o,\rho) \in \L^2 \g^* \times \L^3 \g^* $ of stable forms which
satisfy
\begin{equation}
  \label{eq}
  \o \wedge \rho = 0, \quad d\o^2 = 0, \quad d\rho=0   
\end{equation}
and which induce a Riemannian metric. Working in a basis with fixed
Lie bracket, which determines the exterior derivative completely in
the left-invariant case, two of the equations are quadratic and one is
linear in the coefficients of $\o$ and $\rho$. For each case
separately in a standard basis, large families of solutions of the
equations \eqref{eq} can be constructed with the help of a computer
algebra system, for instance Maple. However, even after Maple was
taught to compute the induced metric, finding a solution inducing a
positive definite metric required a certain persistence, in particular
for the non-unimodular direct sums. We remark that in each case, all
solutions of \eqref{eq} in a small neighbourhood of the constructed
example give rise to a, in most cases rather large, family of
half-flat $\SU(3)$-structures since the condition that the metric is
positive definite is open.

The stable form formalism in dimension six is due to Hitchin,
\cite{Hi1}, \cite{Hi2}, and is explained in section \ref{stable}. The
formalism suggests to consider also half-flat $\SU(p,q)$-structures,
$p+q=3$, with pseudo-Riemannian metrics or even half-flat
$\SL(3,\bR)$-structures where the almost complex structure is replaced
by an almost para-complex structure which is involutive instead of
anti-involutive. In fact, all these structures are described by a pair
of stable forms satisfying \eqref{eq}. An analogue of the Hitchin flow
relates such structures with indefinite metrics to $\G_2^*$-structures
which is elaborated in \cite{CLSS}. More details are recalled in the
preliminary section \ref{prelim}.

In section \ref{su21}, we give an obstruction to the existence of
half-flat $\SU(p,q)$-structures for arbitrary signature which is
stronger than the obstruction established before and applies to 15
classes. Apart from giving an example of a Lie algebra admitting a
half-flat $\SU(1,2)$-structure, but no half-flat $\SU(3)$-structure,
we abstain from completing the classification in the indefinite case
since it would involve constructing approximately 62=78-15-1 explicit
examples of half-flat $\SU(1,2)$-structures.

In section \ref{sl3}, we turn to the para-complex case of
$\SL(3,\bR)$-structures. Again, we give an example of a Lie algebra
admitting a half-flat $\SL(3,\bR)$-structure, but no half-flat
$\SU(p,q)$-structure for any signature. Furthermore, we consider
half-flat $\SL(3,\bR)$-structures on direct sums such that the
summands are mutually orthogonal, as before, and with the additional
assumption, that the metric restricted to each summand is definite. It
turns out that the proof of the classification of ``orthogonal''
half-flat $\SU(3)$-structures in section \ref{orth} generalises with
some sign modifications and we end up with the same list of 15 Lie
algebras. Finally, we consider half-flat $\SL(3,\bR)$-structure such
that both summands are isotropic. The straightforward result is that
such a structure is admitted on a direct sum of three-dimensional Lie
algebras if and only if both summands are unimodular.

The author likes to thank Vicente Cort\'es for suggesting this
project and for many useful discussions.

\section{Preliminaries}
\label{prelim}

\subsection{$\SU(p,q)$-structures and $\SL(m,\bR)$-structures}
\label{Str}
Let $M=M^{2m}$ be an even-dimensional manifold.  An \emph{almost
pseudo-Hermitian structure} $(g,J,\o)$ on $M$ consists of a
pseudo-Riemannian metric $g$, an orthogonal almost complex structure
$J$ and a two-form $\omega= g(.\,, J \,.)$, called the fundamental
two-form. In order to prevent confusion, we point out that many
authors define the fundamental two-form with the opposite sign. An
almost pseudo-Hermitian structure is equivalent to the reduction of
the frame bundle of $M$ to $\U(p,q)$, $p+q=m$, where $(2p,2q)$
is the signature of the metric $g$. A further reduction to
$\SU(p,q)$, i.e.\ an \emph{\SU(p,q)-structure}, is given by a non-trivial
complex $(m,0)$-form $\Psi$ of constant length.

Locally, there is always a pseudo-orthonormal frame $\{ \eta_1, \dots
, \eta_{2m}\}$ such that $J\eta_i = \eta_{i+m}$ and
$\sigma_i=g(\eta_i,\eta_i)=\pm 1$ for $i=1, \dots, m$ and
\[ \o = - \sum_{i=1}^m \sigma_i \, \eta^{i(i+m)} \quad, \qquad \Psi =
(\eta^1 + \i\, \eta^{1+m}) \^ \dots \^ (\eta^m + \i\, \eta^{2m}),\]
where upper indices denote dual basis vectors and $\eta^{ij}$ stands
for the wedge product of $\eta^i$ and $\eta^j$.

Similarly, an \emph{almost para-Hermitian structure} $(g,J,\o)$ on $M$
consists of a neutral metric $g$, an anti-orthogonal para-complex
structure $J$ and the fundamental two-form $\omega = g(\, . \,, J
.\,)$. We recall that an almost para-complex structure $J$ is a
section in $\mathrm{End}(TM)$ such that $J^2 = {\rm id}_{TM}$ and the $\pm
1$-eigenbundles $TM^\pm$ with respect to $J$ have dimension $m$. An
almost para-Hermitian structure is equivalent to a
$\GL(m,\bR)$-structure where $\GL(m,\bR)$ acts reducibly on $T_p M =
T_p M^+ \oplus T_p M^-$ for all $p\in M$. A recent survey on
para-complex geometry is for instance contained in \cite{AMT}.

We denote by $C$ the para-complex numbers $a + {\rm e} b, \: {\rm e}^2
= 1, \: a, b \in \bR$, and by $\Omega^{k,l} M$ the bi-grading induced
by the decomposition of the para-complexification $TM \otimes C$ into
the $\pm {\rm e}$-eigenspaces of $J$. In analogy to the almost
pseudo-Hermitian case, an $(m,0)$-form $\Psi$ of constant
\emph{non-zero} length defines a reduction of the structure group from
$\GL(m,\bR)$ to $\SL(m,\bR)$.

Stressing the similarity to the almost pseudo-Hermitian situation, we
can choose a local pseudo-orthonormal frame $\{ \eta_1, \dots ,
\eta_{2m}\}$ such that $J\eta_i = \eta_{i+m}$ and $g(\eta_i,\eta_i) =
- g(\eta_{i+m},\eta_{i+m})= 1$ for $i=1, \dots, m$ and moreover,
\[ \o = - \sum_{i=1}^m \eta^{i} \wedge \eta^{i+m} \quad, \qquad \Psi =
(\eta^1 + \ep \, \eta^{1+m}) \^ \dots \^ (\eta^m + \ep \,
\eta^{2m}).\] Alternatively, a local frame $\{ \xi_1, \dots ,
\xi_{2m}\}$ can always be chosen such that
\begin{eqnarray}
  g &=& 2 \sum_{i=1}^m \xi^i \cdot \xi^{i+m} , \qquad J\xi_i = \xi_i, \: J\xi_{i+m}=-\xi_{i+m} \quad \mbox{for $i=1, \dots, m$,} \nonumber\\
  \o &=& - \sum_{i=1}^m \xi^{i} \wedge \xi^{i+m}, \quad \Psi = \sqrt 2 \,\{\, (\xi^{1 \dots m} + \xi^{(m+1) \dots 2m}) + \ep (\xi^{1 \dots m} - \xi^{(m+1) \dots 2m}) \,\}\,.  \label{paranormal}
\end{eqnarray}
We will need the following formula, which is easily verified in the
given local frames.
\begin{lemma}
  On an almost pseudo-Hermitian or almost para-Hermitian manifold
  $(M^{2m},g,J,\o)$, the identity
  \begin{eqnarray}
    \label{gonL1}
    \a \wedge J^* \b \wedge \o^{m-1} &=& \frac{1}{m} g(\a,\b) \o^m 
  \end{eqnarray}
  holds for all one-forms $\a$, $\b$.
\end{lemma}

\subsection{Stable forms in dimension six}
\label{stable}
A $p$-form on a vector space is called \emph{stable} if its orbit under
$\GL(V)$ is open \cite{Hi1}. We will frequently use the properties of
stable forms in dimension six and recall the basic facts omitting the
proofs which can be found in \cite{Hi2} and \cite{CLSS}.

Let $V$ be a six-dimensional oriented vector space and $\kappa$ the
canonical isomorphism
\begin{eqnarray*}
  \kappa \, :\, \L^5 V^* \rightarrow V \ot \L^6 V^* , \: \xi \mapsto X \ot \nu \quad \mbox{with  } X \hook \nu = \xi.
\end{eqnarray*}
For every three-form $\rho \in \L^3 V^*$, one can define
\begin{eqnarray}
  \label{K}
  K_{\rho}(v) =& \kappa((v \hook \rho) \wedge \rho) & \in V \ot \L^6 V^*, \\
  \label{la}
  \lambda(\rho) =& \frac{1}{6} {\rm tr}\, K_{\rho}^2 & \in (\L^6 V^*)^{\ot 2} ,\\
  \label{phi}
  \phi(\rho) =& \sqrt{\lambda(\rho)} & \in \L^6 V^*,
\end{eqnarray}
where the positively oriented square root is chosen.  In fact, the
three-form $\rho$ is stable if and only if $\lambda(\rho) \ne 0$. For
a stable three-form $\rho$, we define
\begin{eqnarray}
  \label{J}
  J_{\rho} =& \frac{1}{\phi(\rho)} K_{\rho} & \in {\rm End}(V),
\end{eqnarray}
which is a complex structure if $\lambda(\rho)<0$ and a para-complex
structure for $\lambda(\rho)>0$. Moreover, the form $\rho + \i
J_\rho^* \rho$, or $\rho + \ep J_\rho^* \rho$, respectively, is a
(3,0)-form with respect to $J_\rho$.
\begin{lemma}
  \label{JonL1lemma}
  The (para-)complex structure $J_\rho$ induced by a stable three-form
  $\rho$ acts on one-forms by the formula
  \begin{eqnarray}
    \label{JonL1}
    J_{\rho}^* \alpha (v) \,  \phi(\rho) &=& \alpha \wedge (v \hook \rho) \wedge \rho, \qquad v\in V, \: \a \in V^*.
  \end{eqnarray}
\end{lemma}
\begin{proof}
  The formula follows directly from the definition since we have
  \begin{eqnarray*}
    \alpha \wedge (v \hook \rho) \wedge \rho 
    \stackrel{\eqref{K}}{=} \a \wedge \kappa^{-1}(K_\rho (v))
    \stackrel{\eqref{J}}{=} \a \wedge ( J_{\rho}(v) \hook \phi(\rho) )
    = \alpha (J_{\rho} v) \, \phi(\rho) = J_{\rho}^* \alpha (v) \,  \phi(\rho)
  \end{eqnarray*}
  for all $v\in V$ and all $\a \in V^*$.
\end{proof}
A two-form $\omega \in L^2 V^*$ in dimension six is stable if and only
if it is non-degenerate, i.e.
\[\phi(\o) = \frac{1}{6} \o^3 \ne 0.\]
A pair $(\o,\rho) \in \L^2 V^* \times \L^3 V^*$ of stable forms is
called \emph{compatible} if
\begin{eqnarray}
  \label{comp1}
  \o \^ \rho = 0 &\iff& \omega(.\,,J_\rho \,. ) = - \omega(J_\rho\,.\, , \, .)
\end{eqnarray}
and \emph{normalised} if
\begin{eqnarray}
  \label{comp2}
  \phi(\rho) = \pm 2 \phi(\o) &\iff& J_\rho^* \rho \wedge \rho = \pm \frac{2}{3}\, \omega^3.
\end{eqnarray}
The choice of the sign $\pm$ in the normalisation condition determines
in particular the orientation which is needed to uniquely define
$\phi(\rho)$ and the induced (para-)complex structure $J_\rho$.  A
compatible and normalised pair induces a pseudo-Euclidean metric
\begin{equation} \label{inducedmetric} g = g_{(\o,\rho)}= \varepsilon
  \, \omega(. \,,J_\rho\,.).
\end{equation}
By compatibility, the induced (para-)complex structure $J_\rho$ is
(anti-)orthogonal with respect to this induced metric and the
stabiliser of a compatible and normalised pair is
\begin{eqnarray*}
  \mathrm{Stab}_{\GL(V)}(\rho,\omega) \cong
  \begin{cases}
    \SU(p,q), \, p+q =3, & \mbox{if $\lambda(\rho) < 0$,} \\
    \SL(3,\mathbb R)\,, & \mbox{if $\lambda(\rho) > 0$.}
  \end{cases}
\end{eqnarray*}
In particular, the conventions are chosen such that $\phi(\rho) = + 2
\phi(\o)$ if the induced metric is positive definite which is in fact
the motivation for the sign convention $\o=g(.\,,J\,.)$.

\subsection{Half-flat structures}
\label{hafla}
Let $M$ be a six-manifold. We call $\SU(p,q)$-structures, $p+q=3$, and
$\SL(3,\bR)$-structures defined by tensors $(g,J,\o,\Psi)$
\emph{normalised} if
\[ \Im \Psi \wedge \Re \Psi = \pm \frac{2}{3}\, \omega^3.\] This can
always be achieved by rescaling the length of $\Psi$ which is constant
and non-zero by definition. In fact, the local frames given in section
\ref{Str} are already normalised. Furthermore, we call a $p$-form
$\rho$ on a manifold stable if $\rho_p$ is stable on $T_pM$ for all
$p\in M$. With this terminology, the discussion of stable forms in
dimension six can be applied to $\SU(p,q)$-structures $(g,J,\o,\Psi)$,
$p+q =3$ and $\SL(3,\bR)$-structures as follows.

\begin{proposition}
  Let $M$ be a six-manifold.
  \begin{enumerate}[(i)]
  \item There is a one-to-one correspondence between normalised
    $\SU(p,q)$-structures $(g,J,\o,\Psi)$, \\ $p+q =3$, on $M$ and
    pairs $(\o,\rho) \in \Omega^2 M \times \Omega^3 M$ of stable forms
    which are everywhere compatible and normalised and satisfy
    $\lambda(\rho_p)<0$ for all $p \in M$.
  \item There is a one-to-one correspondence between normalised
    $\SL(3,\bR)$-structures $(g,J,\o,\Psi)$ on $M$ and pairs
    $(\o,\rho) \in \Omega^2 M \times \Omega^3 M$ of stable forms which
    are everywhere compatible and normalised and satisfy
    $\lambda(\rho_p)>0$ for all $p \in M$.
  \end{enumerate}
\end{proposition}

An $\SU(p,q)$-structure, $p+q=3$, or an $\SL(3,\bR)$-structure,
defined by a pair of forms $(\o,\rho)$, is called \emph{half-flat} if
\begin{equation}
  \label{hf}
  d\rho =0 \, , \qquad d \o^2 = 0. 
\end{equation}
We will speak of a \emph{half-flat structure} if the sign of $J^2$ and
the signature of $g$ are not important.

Moreover, left-invariant half-flat structures $(\omega,\rho)$ on Lie
groups $G$ are in one-to-one correspondence with pairs $(\omega,\rho)
\in \L^2\g^* \times \L^3 \g^*$ of stable forms on the corresponding
Lie algebra $\g$ satisfying the exterior system \eqref{hf} and $\o \^
\rho =0$. Therefore, we denote a pair $(\omega,\rho) \in \L^2\g^*
\times \L^3 \g^*$ with these properties as a \emph{half-flat structure
  on a Lie algebra} and the existence problem of left-invariant
half-flat structures on Lie groups reduces to the existence of
half-flat structures on Lie algebras.

\subsection{Three-dimensional Lie algebras}
\label{3d}
Let $\g$ be the Lie algebra of an $n$-dimensional real Lie group
$G$. Identifying $\g$ with the Lie algebra of left-invariant vector
fields on $G$, the formula
\begin{equation*}
  d \a (X,Y) = - \a([X,Y]), \quad \a \in \g^* \, , \; X,Y \in \g,
\end{equation*}
shows that the exterior derivative of $G$ restricted to left-invariant
one-forms contains the same information as the Lie bracket. Since the
Jacobi identity is equivalent to $d^2=0$, we have a complex $(\L^*
\g^*,d)$. Its cohomology $H^*(\g)$ is the Chevalley-Eilenberg or Lie
algebra cohomology for the trivial representation.

Recall that a Lie algebra $\g$ is called \emph{unimodular} if the
trace of the adjoint representation $\rm{ad}_{X}$ vanishes for all $X
\in \g$.
\begin{lemma}
  \label{uni_char}
  The following conditions are equivalent for an $n$-dimensional Lie
  algebra.
  \begin{enumerate}[(i)]
  \item $\g$ is unimodular
  \item All $(n-1)$-forms on $\g$ are closed.
  \item $H^n(\g) = \bR$
  \item Let $\{ c_{ij}^k \}$ denote the structure constants with
    respect to a basis $\{ \e^i \}$ of $\g^*$ which are defined by $d
    \e^k = \sum_{i<j} c_{ij}^k \e^{ij}$. Then, it holds $\sum_{k=1}^n
    c_{k,m}^k = 0$ for $1 \le m \le n$.
  \item The associated connected Lie groups $G$ are unimodular, i.e.\
    the Haar measure of $G$ is bi-invariant.
  \end{enumerate}
\end{lemma}
Unimodularity is a necessary condition for the existence of a
co-compact lattice, see for instance \cite{Mi}, in dimension three it
is also sufficient. Indeed, the closed three-manifolds of the form
$\Gamma \backslash G$ where $G$ is a Lie group with lattice $\Gamma$
are classified in \cite{RV}. Since a direct sum $\g_1 \oplus {\g_2}$
of Lie algebras is unimodular if and only if both $\g_1$ and ${\g_2}$
are unimodular, a direct product $G_1 \times G_2$ of three-dimensional
Lie groups admits a co-compact lattice if and only if it is
unimodular.
\begin{lemma} \label{omegasquare} Let $\g_1 \oplus {\g_2}$ be the
  direct sum of two Lie algebras of dimension three. Moreover, let
  $\omega$ be a non-degenerate two-form in $\Lambda^2 (\g_1 \oplus
  {\g_2})^* = \Lambda^2 \g_1^* \oplus (\g_1^* \otimes \g_2^*)
  \oplus \Lambda^2 \g_2^*$ such that the projections of $\omega$ on
  $\Lambda^2 \g_1^*$ and $\Lambda^2 \g_2^*$ vanish.  Then
  $\omega^2$ is closed if and only if both $\g_1$ and ${\g_2}$ are
  unimodular.
\end{lemma}
\begin{proof}
  Since $\omega \in \g_1^* \otimes \g_2^*$ is non-degenerate, we
  can always choose bases $\{ \e^i \}$ of $\g_1^*$ and $\{ \f^i \}$
  of $\g_2^*$ such that $\omega = \sum_{j=1}^3 \e^j
  \f^j$. Therefore, we have
  \[ \omega^2 = - 2 \sum_{i<j} \e^{ij} \f^{ij} \quad \Rightarrow \quad
  -\frac{1}{2} \, d \omega^2 = \sum_{i<j} d\, (\e^{ij}) \^ \f^{ij} +
  \sum_{i<j} \e^{ij} \^ d\, (\f^{ij}). \] By Lemma \ref{uni_char},
  both $\g_1$ and ${\g_2}$ are unimodular if and only if all two-forms
  $\e^{ij}$ and $\f^{ij}$ are closed. Since the sum is a direct sum of
  Lie algebras, the assertion follows immediately.
\end{proof}
In the following chapter, we need to determine in which isomorphism
class a given three-dimensional Lie algebra lies. All information we
need, including proofs, can be found in \cite{Mi}. We summarise the
results in two propositions. Recall that a Euclidean cross product in
dimension three is determined by a scalar product and an orientation.
\begin{proposition}[Unimodular case] \label{uni3d} Let $\g$ be a
  three-dimensional Lie algebra and choose a scalar product and an
  orientation.
  \begin{enumerate}[(a)]
  \item There is a uniquely defined endomorphism $L$ of $\g$ such that
    $[u,v] = L (u \times v)$.
  \item The Lie algebra $\g$ is unimodular if and only if $L$ is
    self-adjoint.
  \item If $\g$ is unimodular, the isomorphism class of $\g$ is
    characterised by the signs of the eigenvalues of $L$. It can be
    achieved that there is at most one negative eigenvalue of $L$ by
    possibly changing the orientation.
  \end{enumerate}
\end{proposition}

\begin{table}[ht]
  \renewcommand{\arraystretch}{1.6}
  \setlength{\tabcolsep}{0.3cm}
  \begin{tabular}{c|c|c|c}
    & Bianchi type & Eigenvalues of L & Standard Lie bracket \\\hline \hline
    $\su(2) \cong \mathfrak{so}(3)$ & IX & (+,+,+) & $d\e^1\! =\! \e^{23}, \: d\e^2\! = \!\e^{31}, \: d\e^3 \!=\! \e^{12}$ \\ \hline
    $\sl(2,\bR) \cong \mathfrak{so}(1,2)$ & VIII & (+,+,-) & $d\e^1 = \e^{23}, \: d\e^2 = \e^{31}, \: d\e^3 = \e^{21}$ \\ \hline
    $\mathfrak e(2)$ & VII$\mbox{}_0$ & (+,+,0)  & $d\e^2 = \e^{31}, \: d\e^3 = \e^{12}$ \\ \hline
    $\mathfrak e(1,1)$ & VI$_0$ & (+,-,0)  & $d\e^2 = \e^{31}, \: d\e^3 = \e^{21}$ \\ \hline
    $\h_3$ & II & (+,0,0)  & $d\e^3 = \e^{12}$ \\ \hline
    $\bR^3$ & I & (0,0,0)  & abelian \\ 
  \end{tabular}\\[2ex]
  \caption{Three-dimensional unimodular Lie algebras}
  \label{unilist}
\end{table}

Recall that the unimodular kernel of a Lie algebra $\g$ is the kernel
of the Lie algebra homomorphism
\[ \g \rightarrow \bR \; ,\: X \mapsto \rm{tr}(\rm{ad}_X). \]
\begin{proposition}[Non-unimodular case] \label{nonuni3d} Let $\g$ be
  a non-unimodular three-dimensional Lie algebra.
  \begin{enumerate}[(a)]
  \item The unimodular kernel $\mathfrak u$ of $\g$ is two-dimensional
    and abelian.
  \item Let $X \in \g$ such that $\rm{tr}(\rm{ad}_X)=2$ and let
    $\tilde L : \mathfrak u \rightarrow \mathfrak u $ be the
    restriction of $\rm{ad_X}$ to the unimodular kernel $\mathfrak
    u$. If $\tilde L$ is not the identity map, the isomorphism class
    of $\g$ is characterised by the determinant $D$ of $\tilde L$.
  \end{enumerate}
\end{proposition}

\begin{table}[ht]
  \renewcommand{\arraystretch}{1.6}
  \setlength{\tabcolsep}{0.3cm}
  \begin{tabular}{c|c|c|c}
    & Bianchi type & D & Standard Lie bracket  \\\hline \hline
    $\solv_2 \oplus \bR$ & III & 0 & $d\e^2 = \e^{21}$  \\ \hline
    $\solv_3$ & IV & 1 (and $\tilde L \ne$ id) & $d\e^2 = \e^{21} \!+ \e^{31}, \, d\e^3 = \e^{31}$  \\ \hline
    $\solv_{3,1}$  & V & 1 (and $\tilde L =$ id) & $d\e^2 = \e^{21}, \, d\e^3 = \e^{31}$   \\ \hline
    $\solv_{3,\mu}$ {\scriptsize $(-1 <  \mu  < 0)$} & VI & $ D=\frac{4 \mu}{(\mu +1)^2} < 0$ & $d\e^2 = \e^{21}, \, d\e^3 = \mu \e^{31}$   \\ \hline
    $\solv_{3,\mu}$ {\scriptsize $(0 <  \mu  < 1)$} & VI & $0 < D=\frac{4 \mu}{(\mu +1)^2} < 1$ & $d\e^2 = \e^{21}, \, d\e^3 = \mu \e^{31}$   \\ \hline
    $\solv'_{3,\mu}$ {\scriptsize $(\mu > 0)$} & VII & $D=1 + \frac{1}{\mu^2} > 1$ & $d\e^2 = \mu \e^{21} \!+ \e^{13}, \, d\e^3 = \e^{21} \!+ \mu \e^{31}$ \\ 
  \end{tabular}\\[2ex]
  \caption{Three-dimensional non-unimodular Lie algebras}
  \label{nonunilist}
\end{table}

We remark that all three-dimensional Lie algebras are solvable except
for $\su(2)$ and $\sl(2,\bR)$ which are simple. The three-dimensional
Heisenberg algebra $\h_3$ represents the only non-abelian nilpotent
isomorphism class. The two Lie algebras $\mathfrak e(2)$ and
$\mathfrak e(1,1)$ correspond to the groups of rigid motions of the
Euclidean plane $\bR^2$ and of the Minkowskian plane $\bR^{1,1}$,
respectively. The names for the non-unimodular Lie algebras are taken
from \cite{GOV} and the Bianchi types are defined in the original
classification by Bianchi from 1898, \cite{B1}, see \cite{B2} for an
English translation.

\section{Classification of direct sums admitting a half-flat
  $\SU(3)$-structure such that the summands are orthogonal}
\label{orth}
A Hermitian structure on a $2m$-dimensional Euclidean vector space
$(V,g)$ is given by an orthogonal complex structure $J$. As before, we
define the fundamental two-form by $\omega = g(. \,, J .\,)$. The
following Lemma is crucial for the proof of the first classification
result.
\begin{lemma}
  \label{u3ortho}
  Let $(V_1,g_1)$ and $(V_2,g_2)$ be three-dimensional Euclidean
  vector spaces and let $(g,J,\omega)$ be a Hermitian structure on the
  orthogonal product $(V_1 \oplus V_2, g = g_1 + g_2)$. There are
  orthonormal bases $\{ e_1, e_2, e_3 \}$ of $V_1$ and $\{ f_1, f_2,
  f_3 \}$ of $V_2$ which can be joined to an orthonormal basis of $V_1
  \oplus V_2$ such that
  \begin{equation}
    \label{omeganormalperp}
    \omega = a \, \e^{12} + \sqrt{1 - a^2} \, \e^{1} \f^1 + \sqrt{1 -a^2} \, \e^2 \f^2 + \e^3 \f^3 - a \, \f^{12}
  \end{equation}
  for a real number a with $-1 < a \le 1$.
\end{lemma}
\begin{proof}
  Let $\{ e_1, e_2, e_3 \}$ and $\{ f_1, f_2, f_3 \}$ be orthonormal
  bases of $V_1$ and $V_2$, respectively. The group $O(3) \times O(3)$ acts
  transitively on the pairs of orthonormal bases. Let $\Omega$ be the
  Gram matrix of the two-form $\omega$ with respect to our
  basis. Writing the upper right block of $\Omega$ as a product of an
  orthogonal and positive semi-definite matrix and acting with an
  appropriate pair of orthogonal matrices, we find an orthonormal
  basis and nine real parameters such that
  \[ \Omega = \left( \begin{matrix}
      0 & y_1 & y_2 & x_1 & 0 & 0 \\
      -y_1 & 0 & y_3 & 0 & x_2 & 0 \\
      -y_2 & -y_3 & 0 & 0 & 0 & x_3\\
      -x_1 & 0 & 0 &  0 & z_1 & z_2 \\
      0 & -x_2 & 0 & -z_1 & 0 & z_3 \\
      0 & 0 & -x_3 & -z_2 & -z_3 & 0
    \end{matrix} \right)
  \]
  with $x_i \ge 0$ for all $i$ and det$(\Omega) \ne 0$.

  Since $\omega =g(\, . \,, J . \,)$, the matrix $\Omega$ with respect
  to an orthonormal basis has to be a complex structure, i.e.\
  $\Omega^2 = -\id$, where $\id$ denotes the identity matrix. In our
  basis, $\Omega^2$ is {\small
    \[ \left(
      \begin{matrix}
        -y_{1}^2-y_{2}^2-x_{1}^2  & -y_{2} y_{3} & y_{1} y_{3} & 0 & y_{1} x_{2}+x_{1} z_{1} & y_{2} x_{3}+x_{1} z_{2} \\
        -y_{2} y_{3} & -y_{1}^2-y_{3}^2-x_{2}^2  & -y_{1} y_{2} & -y_{1} x_{1}-x_{2} z_{1} & 0 & y_{3} x_{3}+x_{2} z_{3} \\
        y_{1} y_{3} & -y_{1} y_{2} & -y_{2}^2-y_{3}^2-x_{3}^2  & -y_{2} x_{1}-x_{3} z_{2} & -y_{3} x_{2}-x_{3} z_{3} & 0 \\
        0 & -y_{1} x_{1}-x_{2} z_{1} & -y_{2} x_{1}-x_{3} z_{2} & -x_{1}^2-z_{1}^2-z_{2}^2  & -z_{2} z_{3} & z_{1} z_{3} \\
        y_{1} x_{2}+x_{1} z_{1} & 0 & -y_{3} x_{2}-x_{3} z_{3} & -z_{2} z_{3} & -x_{2}^2-z_{1}^2-z_{3}^2   & -z_{1} z_{2} \\
        y_{2} x_{3}+x_{1} z_{2} & y_{3} x_{3}+x_{2} z_{3} & 0 & z_{1}
        z_{3} & -z_{1} z_{2} & -x_{3}^2-z_{2}^2-z_{3}^2
      \end{matrix} \right).
    \]
  } We end up with a set of 18 quadratic equations (and one
  inequality) and determine all solutions modulo the action of $O(3)
  \times O(3)$ and an exchange of the summands.

  On the one hand, assume $y_i = 0$ for all $i$. It follows that all
  equations are satisfied if and only if $x_i = 1$ and $z_i = 0$ for
  all $i$. In this case, the two-form $\omega$ is in the normal form
  $(\ref{omeganormalperp})$ with $a=0$.

  On the other hand, assume that one of the $y_i$ is different from
  zero, say $a:=y_1 \ne 0$ without loss of generality. Inspecting the
  first two terms of the third line of $\Omega^2$, we observe $y_2 =
  y_3 = 0$. Since $x_i \ge 0$, the first three elements on the
  diagonal enforce $x_3=1$, $x_1=x_2=\sqrt{- a^2 + 1}$ and $|a| \le
  1$. But $x_3=1$ and $y_2=y_3=0$ imply that $z_2=z_3=0$ due to row 3,
  terms 4 and 5. If $|a|<1$ and thus $x_1=x_2 > 0$, the term in row 1
  and column 5 enforces $z_1 = -a$. Obviously, all equations are
  satisfied and $\omega$ is in the normal form
  $(\ref{omeganormalperp})$. Finally, if $|a|=1$, we have immediately
  $x_1=x_2=0$ and $|z_1|=1$ and all equations are satisfied
  again. Since changing the signs of the base vectors $e_1$ and $f_1$
  is an orthogonal transformation which does not change $x_3$, we can
  obtain the normal form $(\ref{omeganormalperp})$ for $a = 1$. Since
  we found all solutions to the 18 equations and the two-form $\omega$
  is non-degenerate for all values of $a$, the Lemma is proven.
\end{proof}
We call the Hermitian structure of type I if it admits a basis with
$a=0$ and of type II if it admits a basis with $a\ne 0$.
\begin{theorem}
  \label{theo_orth}
  Let $\g=\g_1 \op \g_2$ be a direct sum of three-dimensional Lie
  algebras $\g_1$ and $\g_2$. 

  The Lie algebra $\g$ admits a half-flat $\SU(3)$-structure
  such that $\g_1$ and $\g_2$ are mutually orthogonal and such that
  the underlying Hermitian structure is of type I if and only if
  \begin{align*}
    \mbox{(i)} & \quad \mbox{$\g_1=\g_2$ and both are unimodular or } \\
    \mbox{(ii)} & \quad \mbox{$\g_1$ is non-abelian unimodular and $\g_2$ abelian or vice versa.} 
  \end{align*}
  Moreover, the Lie algebra $\g$ admits a half-flat $\SU(3)$-structure
  such that $\g_1$ and $\g_2$ are mutually orthogonal and such that
  the underlying Hermitian structure is of type II if and only if the
  pair $(\g_1,\g_2)$ or $(\g_2,\g_1)$ is contained in the following
  list:
  \begin{align*}
    &&& (\mathfrak e(1,1),\mathfrak e(1,1)), \\
    &&& (\mathfrak e(2),\bR \oplus \solv_{2}), \\
    &&& (\su(2),\solv_{3,\mu})  && \mbox{for $0 < \mu \le 1$,} && \\
    &&& (\sl(2,\bR),\solv_{3,\mu}) && \mbox{for $-1 < \mu < 0$.} &&
  \end{align*}
\end{theorem}

\begin{proof}
  Given an arbitrary (almost) Hermitian structure $(g, J, \omega )$
  such that $\g_1$ and $\g_2$ are orthogonal, we can use Lemma
  \ref{u3ortho} and choose an orthonormal basis $\{ e_1, e_2, e_3,
  f_1, f_2, f_3 \}$ of $\g_1 \oplus \g_2$ such that $\{ e_1, e_2, e_3
  \}$ spans $\g_1$, $\{ f_1, f_2, f_3 \}$ spans $\g_2$ and
  \begin{equation}
    \omega = a \, \e^{12} + \sqrt{1 - a^2} \, \e^{1} \f^1 + \sqrt{1 -a^2} \, \e^2 \f^2 + \e^3 \f^3 - a \, \f^{12}
  \end{equation}
  for a real number a with $-1 < a \le 1$. The reductions from $\U(3)$
  to $\SU(3)$ are parameterised by the space of complex-valued
  $(3,0)$-forms $\Psi=\psi+\i \phi$ which is complex
  one-dimensional. We remark that, working on a vector space, the
  length normalisation of the $(3,0)$-form is not important for the
  existence question. The Lie bracket of the direct sum $\g_1 \oplus
  \g_2$ is encoded in the 18 structure constants of $\g_1$ and $\g_2$:
  \[ d \e^i= c_{j,k}^{i} \e^{jk} \quad \mbox{and} \quad d \f^i=
  c_{j+3,k+3}^{i+3} \f^{jk} \quad \mbox{with i,j,k} \in \{ 1,2,3
  \}. \] Therefore, our ansatz includes 21 parameters consisting of 18
  structure constants, two real parameters defining an arbitrary
  $\SU(3)$ reduction and the parameter $a$. Our strategy is to find
  all solutions of the equations defining half-flatness
  \begin{eqnarray*}
    d \omega^2 = 0  & \mbox{and} & d \psi  = 0
  \end{eqnarray*}
  and the Jacobi identity $d^2 = 0$ and to determine the isomorphism
  classes of the solutions if necessary.

  Type I: Assume first that $a=0$. Due to Lemma (\ref{omegasquare}),
  the first half-flat equation $d \omega^2=0$ is satisfied if and only
  if both $\g_1$ and $\g_2$ are unimodular. It remains to solve the
  second half-flat equation for unimodular summands. Since we have
  $J(f_i)=e_i$ in our basis for $a=0$, the dual vectors satisfy $\e^i
  \circ J = \f^i$. Therefore, the complex-valued form
  \begin{eqnarray*}
    \Psi_0 & = & \psi_0 + \i \phi_0 = (\e^1 - \i \e^1 \circ J) \wedge (\e^2 - \i \e^2 \circ J) \wedge (\e^3 - \i \e^3 \circ J) \\
    & = & \e^{123} - \e^{1}\f^{23} - \e^{2}\f^{31} - \e^{3}\f^{12} + \i (\f^{123} - \e^{12}\f^{3} - \e^{31}\f^{2} - \e^{23}\f^{1})
  \end{eqnarray*}
  is a $(3,0)$-form with respect to $J$.  By multiplying $\Psi_0$ with
  a non-zero complex number $\xi_1 + \i \xi_2$, we obtain all
  $(3,0)$-forms. Their real part is $\psi = \xi_1 \psi_0 - \xi_2
  \phi_0$. Considering that all two-forms on both $\g_1$ and $\g_2$
  are closed, we compute the exterior derivative of $\psi$:
  \begin{eqnarray*}
    d \psi &=&
    -(\xi_1 c_{1,2}^{1}-\xi_2 c_{5,6}^{6}) \:\e^{12}\f^{23}
    -(\xi_1 c_{2,3}^{1}-\xi_2 c_{5,6}^{4}) \:\e^{23}\f^{23}
    -(\xi_1 c_{3,1}^{1}-\xi_2 c_{5,6}^{5}) \:\e^{31}\f^{23}\\
    &&-(\xi_1 c_{1,2}^{2}-\xi_2 c_{6,4}^{6}) \:\e^{12}\f^{31}
    -(\xi_1 c_{2,3}^{2}-\xi_2 c_{6,4}^{4}) \:\e^{23}\f^{31}
    -(\xi_1 c_{3,1}^{2}-\xi_2 c_{6,4}^{5}) \:\e^{31}\f^{31}\\
    &&-(\xi_1 c_{1,2}^{3}-\xi_2 c_{4,5}^{6}) \:\e^{12}\f^{12}
    -(\xi_1 c_{2,3}^{3}-\xi_2 c_{4,5}^{4}) \:\e^{23}\f^{12}
    -(\xi_1 c_{3,1}^{3}-\xi_2 c_{4,5}^{5}) \:\e^{31}\f^{12}.
  \end{eqnarray*}
  If $\xi_1$ or $\xi_2$ is zero we have obviously $d \psi=0$ if and
  only if one of the summands is abelian. By Lemma (\ref{uni_char}),
  the unimodularity of $\g_2$ is equivalent to $c_{6,4}^6 =
  -c_{5,4}^{5}$, $c_{6,5}^6 = -c_{4,5}^{4}$ and $c_{5,4}^5 =
  -c_{6,4}^{6}$. Therefore, if both $\xi_1$ and $\xi_2$ are different
  from zero, $d \psi$ vanishes if and only if the structure constants
  of $\g_1$ and $\g_2$ coincide up to the scalar $\frac{\xi_1}{\xi_2}$
  and therefore $\g_1 = \g_2$. This comprises all solutions under the
  assumption $a=0$.

  Type II: Assume now that the $\U(3)$-structure satisfies $a \ne
  0$. To improve readability, the abbreviation $b:=\sqrt{1 - a^2}$ is
  introduced.

  With this notation, we compute
  \begin{align*}
    \frac{1}{2} \omega^2 = a \: \e^{123}\f^{3}- a \: \e^{3}\f^{123} -
    \: \e^{12}\f^{12} - b \: \e^{13}\f^{13}- b \: \e^{23}\f^{23} &&&
  \end{align*}
  \vspace{-0.5cm}
  \begin{align*}
    &&&&&&&&\frac{1}{2} d (\omega^2)
    &=(c_{4,6}^{4}+ c_{5,6}^{5} - a c_{1,2}^{3})&\e^{12}\f^{123}\:
    &+\:( b c_{5,6}^{6} - b c_{4,5}^{4} - a c_{1,3}^{3})&\e^{13}\f^{123}&&&&&&&&&\\
    &&&&&&&& &-( a c_{2,3}^{3}+ b c_{4,5}^{5}+ b c_{4,6}^{6})&\e^{23}\f^{123}\:
    &+\:( c_{1,3}^{1}+ c_{2,3}^{2}- a c_{4,5}^{6})&\e^{123}\f^{12}&&&&&&&&&\\
    &&&&&&&& &-( b c_{1,2}^{1}- b c_{2,3}^{3}+ a c_{4,6}^{6})&\e^{123}\f^{13}\:
    &-\:( a c_{5,6}^{6}+ b c_{1,2}^{2}+ b c_{1,3}^{3})&\e^{123}\f^{23}.&&&&&&&&&
  \end{align*}
  We reduce our ansatz to the space of solutions of $d \omega^2=0$ by
  substituting
  \begin{eqnarray*}
    c_{2,3}^{2} &=& a c_{4,5}^{6}-c_{1,3}^{1} \quad , \quad
    c_{2,3}^{3} = b^2 c_{1,2}^{1}- a b c_{4,5}^{5} \quad , \quad
    c_{1,3}^{3} = - b^2 c_{1,2}^{2} - a b c_{4,5}^{4} \quad , \quad\\
    c_{5,6}^{5} &=& a c_{1,2}^{3}-c_{4,6}^{4} \quad , \quad  
    c_{5,6}^{6} = b^2 c_{4,5}^{4} - a b c_{1,2}^{2} \quad \mbox{and} \quad  
    c_{4,6}^{6} = - b^2 c_{4,5}^{5} - a b c_{1,2}^{1} .
  \end{eqnarray*}
  In our basis, we have $\e^1 \circ J = b \f^1 + a \e^2$, $\e^3 \circ
  J = \f^3$ and $\f^2 \circ J = - b \e^2 + a \f^1$. Using this, we
  compute a (3,0)-form $ \Psi_0 = \psi_0 + \i \phi_0 $ with
  \begin{eqnarray*}
    \psi_0 & =&  + b\:\f^{123} - b \: \e^{12}\f^{3} + \e^{13}\f^{2} - \:\e^{23}\f^{1} + a \:\e^{1}\f^{13} + a \:\e^{2}\f^{23}  \\
    \phi_0 &=& - b\:\e^{123} + b \:\e^{3}\f^{12}  - \:\e^{2}\f^{13} +\:\e^{1}\f^{23} - a \:\e^{13}\f^{1} -a \:\e^{23}\f^{2} 
  \end{eqnarray*}
  In the following, we work with the real part $\psi = \xi_1 \psi_0 -
  \xi_2 \phi_0$ of an arbitrary (3,0)-form. By possibly changing the
  roles of the two summands, we can assume that $\xi_1$ is non-zero
  and we normalise our (3,0)-form such that $\xi_1=1$.  The exterior
  derivative of $\psi$ is, after inserting the above substitutions,
  \begin{align*}
    d \psi = a b \, c_{4,5}^{6} \,\e^{123}\f^{3} \, &&- \xi_2 a b \,
    c_{1,2}^{3} \,\e^{3}\f^{123} && - b \,(\,c_{4,5}^{6} + \xi_2
    \,c_{1,2}^{3} )\,\e^{12}\f^{12} &&&&&&&&&&&&&&&&&&&
  \end{align*}
  \vspace{-0.5cm}
  \begin{align*}
    &+(& -\xi_2 a b \, c_{1,2}^{1} && -a^2 b \,c_{1,2}^{2} && -a^3
    c_{4,5}^{4} && +\xi_2 a^2 c_{4,5}^{5} &) \,\e^{1}\f^{123} &&&
    \\ 
    &+(& a^2 b \,c_{1,2}^{1} && -\xi_2 a b \,c_{1,2}^{2} && -\xi_2 a^2
    c_{4,5}^{4} &&-a^3 c_{4,5}^{5} &) \,\e^{2}\f^{123} &&& \\ 
    &+(& \xi_2 a^3 c_{1,2}^{1} && - a^2 c_{1,2}^{2} &&+ a b \,
    c_{4,5}^{4} &&+\xi_2 a^2 b \,c_{4,5}^{5} &
    )\,\e^{123}\f^{1}&&&\\ 
    &+(&a^2 c_{1,2}^{1} && +\xi_2 a^3 c_{1,2}^{2} &&-\xi_2 a^2 b
    \,c_{4,5}^{4} &&+ a b \, c_{4,5}^{5} &)\,\e^{123}\f^{2} &&&\\
    &+(&a(2 -a^2) \,c_{1,2}^{1} &&+ \xi_2 \, c_{1,2}^{2} && &&+ b^3
    \,c_{4,5}^{5} &) \,\e^{12}\f^{13}&&&\\ 
    &+(&- \xi_2 \, c_{1,2}^{1} &&+ \,a(2 -a^2) \, c_{1,2}^{2} &&- b^3
    \,c_{4,5}^{4} && &) \,\e^{12}\f^{23}&&&\\ 
    &+(& && +\xi_2 b^3 \,c_{1,2}^{2} && \xi_2 a(2 - a^2)\,c_{4,5}^{4}
    &&+c_{4,5}^{5} &)\,\e^{13}\f^{12}\\ 
    &+(& - \xi_2 b^3 \, c_{1,2}^{1} && &&- c_{4,5}^{4} &&+\xi_2 a (2-
    a^2) \, c_{4,5}^{5} &) \,\e^{23}\f^{12}&&&
  \end{align*}
  \vspace{-0.5cm}
  \begin{align*}
    &&&+(& && a\, c_{1,3}^{1} &&+\xi_2 \, c_{1,3}^{2} && &&+\xi_2 a \,
    c_{4,6}^{4} &&+c_{4,6}^{5} && && &) \,\e^{13}\f^{13}\\ 
    &&&+(& \!\!\!\!\!\! \xi_2 a^2 c_{1,2}^{3} && -a \, c_{1,3}^{1} &&
    &&-\xi_2 c_{2,3}^{1}&& -\xi_2 a \,c_{4,6}^{4} && && -c_{5,6}^{4}
    && +a^2
    c_{4,5}^{6} &) \,\e^{23}\f^{23} \\
    &&&+(& a \, c_{1,2}^{3} &&-\xi_2 \, c_{1,3}^{1} &&+a \,
    c_{1,3}^{2} && && -c_{4,6}^{4} && && + \xi_2 a \, c_{5,6}^{4} &&
    &) \,\e^{13}\f^{23}\\ 
    &&&+(& && -\xi_2 \, c_{1,3}^{1} && && +a \,c_{2,3}^{1} &&
    -c_{4,6}^{4} && +\xi_2 a \, c_{4,6}^{5} && && +\xi_2 a \,
    c_{4,5}^{6} &) \,\e^{23}\f^{13}. \\
  \end{align*}
  We need to determine all solutions of the coefficient equations of
  $d\psi=0$. First of all, we observe that the variables
  $c_{1,2}^{1}$, $c_{1,2}^{2}$, $c_{4,5}^{4}$ and $c_{4,5}^{5}$ are
  subject to eight linear equations and claim that there is no
  non-trivial solution of this linear system. Indeed, the determinant
  of the four by four coefficient matrix of the first four equations
  is $a^4 (a^2 \xi_2^2+1) (a^2+\xi_2^2) (a^2+b^2)^2 = a^4 (a^2
  \xi_2^2+1) (a^2+\xi_2^2)$ and thus never vanishes for $a \ne 0$. To
  deal with the remaining eight structure constants, subject to seven
  equations, we treat three cases separately.
  \begin{enumerate}[(a)]
  \item Assume first that $b \ne 0$ and $\xi_2 \ne 0$, i.e.\ $0< |a| <
    1$. Obviously, we have $c_{1,2}^{3}=0$ and $c_{4,5}^{6}=0$ by the
    vanishing of the first three coefficients. Moreover, applying easy
    row transformations to the remaining four equations, we observe
    that it holds necessarily $c_{1,3}^{2} = c_{2,3}^{1}$ and
    $c_{4,6}^{5} = c_{5,6}^{4}$. Considering this, $d\psi=0$ is
    finally satisfied if and only if
    \begin{eqnarray*}
      &s := c_{4,6}^{4} =  \frac{a\, (\xi_2^2-1) \, c_{5,6}^{4} - (a^2 + \xi_2^2) \,c_{1,3}^{1}}{\xi_2 (a^2 +1)}  \, , 
      \quad t := c_{2,3}^{1} =  - \frac{(\xi_2^2 a^2 + 1) \, c_{5,6}^{4} + a (1-\xi_2^2) \, c_{1,3}^{1}}{\xi_2 (a^2+1)} .&
    \end{eqnarray*}
    Applying all substitutions, the set of solutions of the two
    half-flat equations is parameterised by the four parameters $a$,
    $\xi_2$,
    \[ p:=c_{1,3}^{1} \quad \mbox{and} \quad q:=c_{5,6}^{4}. \] In
    order to determine the isomorphism class of $\g_1$ and $\g_2$ for
    all solutions, we apply Propositions \ref{uni3d} and
    \ref{nonuni3d}.  We choose orientations on $\g_1$ and ${\g_2}$
    such that $e_1 \times e_2 = - e_3$ and $e_4 \times e_5 = -
    e_6$. Let $L_{\g_1}$ and $L_{\g_2}$ denote the matrices
    representing the endomorphisms defined in Proposition \ref{uni3d}
    with respect to our bases. On the set of solutions, they simplify
    to
    \[ L_{\g_1} = \mmm{t}{-p}{0}{-p}{-t}{0}{0}{0}{0} \quad \mbox{and}
    \quad L_{\g_2} = \mmm{q}{-s}{0}{-s}{-q}{0}{0}{0}{0}. \] The Jacobi
    identity is already satisfied. Both $L_{\g_1}$ and $L_{\g_2}$ are
    symmetric and in consequence, both summands are unimodular.  The
    eigenvalues of $L_{\g_1}$ and $L_{\g_2}$ are $\{ 0,\pm
    \sqrt{p^2+t^2}\}$ and $\{ 0,\pm \sqrt{s^2+q^2}\}$. Hence, if $p
    \ne 0$ or $q \ne 0$, the Lie algebra $\g_1 \oplus {\g_2}$ is
    isomorphic to $\mathfrak e(1,1) \oplus \mathfrak e(1,1)$ with two
    remaining parameters $\xi_2 \ne 0$ and $0 < |a| < 1$. If $p = 0$
    and $q = 0$, the Lie algebra is abelian.

  \item Now assume $b \ne 0$ and $\xi_2=0$. In this case, the
    equations simplify considerably and the only solution of $d
    \psi=0$ is given by
    \begin{eqnarray*}
      & & c_{4,5}^{6} = 0 \: , \quad c_{4,6}^{4} = a c_{2,3}^{1} \: , \quad  c_{4,6}^{5} = -a c_{1,3}^{1} \: , \quad c_{5,6}^{4} = -a c_{1,3}^{1}  \: , \quad c_{1,2}^{3} = -c_{1,3}^{2}+c_{2,3}^{1}.
    \end{eqnarray*}
    As before, we rename the remaining parameters
    \[ p:=c_{1,3}^{2} \quad \mbox{,} \quad q:=c_{2,3}^{1} \quad
    \mbox{and} \quad r:=c_{1,3}^{1}, \] and have a closer look at
    \[ L_{\g_1} = \mmm{q}{-r}{0}{-r}{-p}{0}{0}{0}{-p+q} \quad
    \mbox{and} \quad L_{\g_2} =
    \mmm{-ar}{-aq}{0}{-ap}{ar}{0}{0}{0}{0}. \] Again, the Jacobi
    identity is already satisfied. The first summand is always
    unimodular, the second summand is unimodular if and only if
    $p=q$. If $p=q$, both matrices are of the same type as in case (a)
    and $\g_1 \oplus {\g_2}$ is isomorphic to $\mathfrak e(1,1) \oplus
    \mathfrak e(1,1)$ or abelian.

    It remains to apply Proposition \ref{nonuni3d} to identify the
    isomorphism class of the solutions with $p \ne q$. Without
    changing the isomorphism class, we can normalise such that
    $p=q+1$. We need to find a vector $X \in {\g_2}$ with
    $\rm{tr}(\rm{ad}_{X})=2$. Since
    $\rm{tr}(\rm{ad}_{f_3})=c_{4,6}^{4} + c_{5,6}^{5} = -a$, we choose
    $X = -\frac{2}{a} f_3$. The unimodular kernel $\mathfrak u$ is
    spanned by $f_1$ and $f_2$ and the restriction of $\rm{ad}_X$ on
    $\mathfrak u$ is represented by the matrix
    \[ \tilde L_{\g_2} = \mm{-2q}{2r}{2r}{2(q+1)} \qquad \mbox{with}
    \quad D=\det(\tilde L_{\g_2})= -4(q(q+1) + r^2) \le 1 . \] If
    $\tilde L_{\g_2}$ is not the identity matrix, the value of $D$
    determines the isomorphism class of ${\g_2}$. However, the
    corresponding class of the unimodular summand $\g_1$ varies with
    the value of $D$. In fact, with $r^2 = -q(q+1)-\frac{1}{4}D$, the
    eigenvalues of $L_{\g_1}$ are $-1$ and $-\frac{1}{2} \pm
    \frac{1}{2} \sqrt{1-D}$. Comparing with the lists in chapter 1, we
    find the remaining classes listed in the theorem.

  \item The last case to be discussed is $b=0$ which corresponds to
    $a=1$.  Now, the equation $d\psi = 0$ is equivalent to
    \begin{eqnarray*}
      && c_{2,3}^{1} = -\xi_2 c_{5,6}^{4}+\xi_2 c_{1,3}^{1}+c_{4,6}^{4} \, , \quad 
      c_{1,2}^{3} = \xi_2 c_{1,3}^{1}+c_{4,6}^{4}-\xi_2 c_{5,6}^{4}-c_{1,3}^{2}, \\
      && c_{4,6}^{5} = -\xi_2 c_{1,3}^{2}-\xi_2 c_{4,6}^{4}-c_{1,3}^{1}  \, , \quad
      c_{4,5}^{6} = \xi_2 c_{4,6}^{4}+\xi_2 c_{1,3}^{2}+c_{5,6}^{4}+c_{1,3}^{1}.
    \end{eqnarray*}
    Considering these substitutions, the Jacobi identity is satisfied
    if and only if
    \begin{equation*}
      \xi_2 c_{4,6}^{4}+\xi_2 c_{1,3}^{2}+c_{5,6}^{4}+c_{1,3}^{1}= 0 \quad \mbox{or} \quad \xi_2 c_{1,3}^{1}+c_{4,6}^{4}-\xi_2 c_{5,6}^{4}-c_{1,3}^{2} = 0.   
    \end{equation*}
    Writing down the matrices $L_{\g_1}$ and $L_{\g_2}$ for both
    cases, it is easy to see that they are of the same form as in case
    (b). Therefore, the possible isomorphism classes of Lie algebras
    are exactly the same as in case (b).
  \end{enumerate}
  Since we have discussed all solutions of the half-flat equations,
  the theorem is proved.
\end{proof}

\section{Classification of direct sums admitting a half-flat $\SU(3)$-structure}

\subsection{Obstructions to the existence of half-flat $\SU(3)$-structures}
\label{sectionobst}
In this section, we establish an obstruction to the existence of
half-flat $\SU(3)$-structures on Lie algebras following the idea of
\cite[Theorem 2]{C}.

We denote by $Z^p$ the space of closed $p$-forms on a Lie algebra and
by $\Ann{W}$ the annihilator of a subspace $W$.
\begin{lemma}
  \label{Jinv}
  Let $\g$ be a six-dimensional Lie algebra and $\g^* = V \op W$ a
  (vector space) decomposition such that $V$ is two-dimensional and
  such that
  \begin{eqnarray}
    \label{*}
    Z^3 &\subset& \L^2 V \^ W \op V \^ \L^2 W.
  \end{eqnarray}
  Then, the subspace $V$ is $J_{\rho}$-invariant for all closed stable
  three-forms $\rho$.
\end{lemma}
\begin{proof}
  Let $\rho \in Z^3$ be stable and $\a \in V$. Since $\dim V = 2$, the
  assumption \eqref{*} implies for all $v \in \Ann{V}$
  \begin{equation*}
    v \hook \rho \in \L^2 V \op V \^ W\, , \qquad \alpha \wedge \rho \in \L^3 V \wedge W \op \L^2 V \^ \L^2 W.
  \end{equation*}
  Therefore, it holds
  \begin{eqnarray*}
    0 = \alpha \wedge (v \hook \rho) \wedge \rho 
    \stackrel{\eqref{JonL1}}{=} J_{\rho}^* \alpha (v) \,  \phi(\rho)
  \end{eqnarray*}
  for all $v \in \Ann{V}$ and, by the definition of the annihilator
  $\Ann{V}$, the subspace $V$ is $J_{\rho}$-invariant.
\end{proof}
\begin{proposition}
  \label{obst}
  Let $\g$ be a six-dimensional Lie algebra and $\g^* = V \op W$ a
  decomposition such that $V$ is two-dimensional and such that
  \begin{eqnarray}
    \label{h03}
    Z^3 &\subset& \L^2 V \^ W \quad \op  V \^ \L^2 W, \\
    \label{h04}
    Z^4 &\subset& \L^2 V \^ \L^2 W  \op  V \^ \L^3 W.
  \end{eqnarray}
  Then, the subspace $V$ is isotropic and $J_\rho$-invariant for every
  half-flat structure $(\o,\rho)$. In particular, the Lie algebra $\g$
  does not admit a half-flat $\SU(3)$-structure.
\end{proposition}
\begin{proof}

  Suppose that $(\omega,\rho)$ is a half-flat structure on $\g$, in
  particular $\rho \in Z^3$ and $\omega^2 \in Z^4$ by definition. By
  Lemma \ref{Jinv}, the subspace $V$ is $J_{\rho}$-invariant. Thus,
  the assumption \eqref{h04} and $\dim V = 2$ imply that
  \begin{eqnarray*}
    0 = \a \^ J_{\rho}^* \b \^ \omega^2 \stackrel{\eqref{gonL1}}{=} \frac{1}{3} \, g(\a,\b) \, \omega^3
  \end{eqnarray*}
  for all $\a,\b \in V$ and $V$ has to be an isotropic subspace of
  $\g^*$. This is of course impossible for definite metrics and there
  cannot exist a half-flat $\SU(3)$-structure.
\end{proof}

\begin{definition}
  Let $\g$ be a Lie algebra. A decomposition $\g^* = V \oplus W$ is
  called a \emph{coherent splitting} if
  \begin{eqnarray}
    \label{coh1}  
    d V &\subset& \Lambda^2 V , \\
    \label{coh2}  
    d W &\subset& \Lambda^2 V \oplus V \wedge W.
  \end{eqnarray}
\end{definition}
\begin{remark}
  The definition can be reformulated into an equivalent dual
  condition:
  \begin{eqnarray*}
    \eqref{coh1} \iff& 0 = d \s (X, .) = - \s([X,.]) \quad \mbox{for all $X \in \Ann{V}, \s \in V$} & \iff [\Ann{V}, \g] \subset \Ann{V},\\
    \eqref{coh2} \iff& 0 = d \s (X, Y) = - \s([X,Y]) \quad \mbox{for all $X,Y \in \Ann{V}, \s \in W$} & \iff [\Ann{V}, \Ann{V}] \subset \Ann{W}.
  \end{eqnarray*}
  In other words, a coherent splitting corresponds to a decomposition
  of $\g$ into an abelian ideal and a vector space complement.
\end{remark}

As elaborated in \cite{C}, a coherent splitting with $\dim V=2$ allows
the introduction of a double complex such that the obstruction
conditions \eqref{h03}, \eqref{h04} can be formulated in terms of the
cohomology of this double complex. However, in the situation we are
interested in, it turns out to be more practical to avoid homological
algebra. Indeed, the verification of the obstruction conditions can be
simplified as follows.  
\begin{lemma}
  \label{splittings}
  Let $\g = \g_1 \oplus \g_2$ be a direct sum of three-dimensional Lie
  algebras.
  \begin{enumerate}[(i)]
  \item Let $\alpha_1 \in \g_1^*$ and $\alpha_2 \in \g_2^*$ be
    one-forms defining $V= \mathrm{span} (\a_1,\a_2)$. Then $\g^* = V \op W$
    is a coherent splitting for any complement $W$ of $V$ if and only
    if the two one-forms $\alpha_i$ are closed and satisfy
    \begin{eqnarray}
      \label{3p3}
      \im (d: \g_i^* \rightarrow \L^2 \g_i^*) &\subset& \alpha_i \wedge \g_i^* \qquad \mbox{for both $i$.}
    \end{eqnarray}
  \item If both summands are non-abelian, \emph{every} coherent
    splitting with $\dim V=2$ is defined by closed one-forms $\alpha_1
    \in \g_1^*$ and $\alpha_2 \in \g_2^*$ satisfying \eqref{3p3}.
  \item There exists a coherent splitting with $\dim V=2$ on $\g$ if
    and only if $\g$ is solvable.
  \item If $\g$ is unimodular, there is no decomposition $\g^* = V \op
    W$ with two-dimensional $V$ satisfying both obstruction conditions
    (\ref{h03}) and (\ref{h04}).
  \end{enumerate}
\end{lemma}
\begin{proof}
  \begin{enumerate}[(i)]
  \item Since both the exterior algebras $\L^* \g_i^*$ are
    $d$-invariant, the condition \eqref{coh1} is satisfied if and only
    if both generators are closed and \eqref{coh2} is equivalent to
    \eqref{3p3}.
  \item Assume that both summands $\g_i$ are not abelian and let a
    coherent splitting be defined by an abelian four-dimensional ideal
    $\Ann{V}$ and a complement. In consequence, both the intersection
    $\Ann{V} \cap \g_i$ and the projection of $\Ann{V}$ on $\g_i$ are
    abelian subalgebras of $\g_i$ for both $i$ and thus at most
    two-dimensional. Since a one-dimensional intersection $\Ann{V}
    \cap \g_i$ would require the projection on the other summand to be
    three-dimensional, it follows that the intersections $\Ann{V} \cap
    \g_i$ have to be exactly two-dimensional. Equivalently, the
    two-dimensional space $V$ is generated by two one-forms $\alpha_1
    \in \g_1^*$ and $\alpha_2 \in \g_2^*$. Now, the assertion follows
    from part (i).
  \item On the one hand, if $\g$ is not solvable, one of the summands
    has to be simple, say $\g_1$. However, the intersection of a
    four-dimensional abelian ideal with $\g_1$ would be zero or
    $\g_1$, both of which is not possible since $\dim \g =6$ and
    since $\g_1$ is not abelian. On the other hand, inspecting the
    list of standard bases in tables \ref{unilist} and
    \ref{nonunilist} reveals that any three-dimensional solvable Lie
    algebra ${\h}$ contains a closed one-form $\alpha$ such that $\im
    \, d \subset \alpha \wedge \h^*$. Therefore, if $\g$ is
    solvable, i.e.\ both summands are solvable, a coherent splitting
    exists by part (i).
  \item Assume that $\g$ is unimodular and let $W$ be an arbitrary
    four-dimensional subspace of $\g^*$. It suffices to show that
    there always exists a closed three-form with non-zero projection
    on $\L^3 W$ \emph{or} a closed four-form with non-zero projection
    on $\L^4 W$. If the projection of $W$ on one of the summands
    $\g_i$ is surjective, every non-zero element of $\L^3 \g_i^*$ is
    closed and has non-zero projection on $\L^3W$. Otherwise, the
    image of the projection of $W$ on either of the summands has to be
    two-dimensional for dimensional reasons. In this case, there is
    always a closed four-form with non-zero projection on $\L^4 W$
    since all four-forms in $\L^2 \g_1^* \^ \L^2 \g_2^*$ are closed by
    unimodularity. This finishes the proof of the lemma.
  \end{enumerate}
\end{proof}

\begin{lemma}
  \label{obst_particular_situation}
  Let $\g = \g_1 \oplus \g_2$ be a direct sum of three-dimensional Lie
  algebras and let $\g^* =V \op W$ be a coherent splitting such that
  $V= \mathrm{span} (\a_1,\a_2)$ is defined by closed one-forms $\a_1
  \in \g_1^*$ and $\a_2 \in \g_2^*$ satisfying \eqref{3p3}. Then,
  the obstruction conditions \eqref{h03} and \eqref{h04} are
  equivalent to the condition that $d$ is injective when restricted to
  $\L^3W$ and $\L^4 W$.
\end{lemma}
\begin{proof}
    The injectivity of $d$ on $\L^3W$ and $\L^4 W$ is obviously
    necessary for \eqref{h03} and \eqref{h04}. With the assumptions,
    it is also sufficient since the coherent splitting satisfies $d W
    \subset V \wedge W$ and $d V =0$ such that the images of $\L^3 W$
    and $\L^4 W$ are linearly independent from the images of the
    complements $\L^2 V \^ W \op V \^ \L^2 W$ and $\L^2 V \^ \L^2 W
    \op V \^ \L^3 W$, respectively.
\end{proof}

\subsection{The classification}
\label{classi}
Using the obstruction established in the previous section, we obtain
the following classification result.
\begin{theorem}
  A direct sum $\g = \g_1 \op \g_2$ of three-dimensional Lie algebras
  admits a half-flat $\SU(3)$-structure if and only if
  \begin{enumerate}[(i)]
  \item $\g$ is unimodular or
  \item $\g$ is not solvable or
  \item $\g$ is isomorphic to $\mathfrak e (2) \op \solv_2 \op \bR$ or
    $\mathfrak e (1,1) \op \solv_2 \op \bR$.
  \end{enumerate}
\end{theorem}
\begin{proof}
  A standard basis of $\g= \g_1 \op \g_2$ will always denote the union
  of a standard basis $\{ \e_1, \e_2, \e_3 \}$ of $\g_1$ and a
  standard basis $\{ \f_1, \f_2, \f_3\}$ of $\g_2$ as defined in
  tables \ref{unilist} and \ref{nonunilist}. For all Lie algebras
  admitting a half-flat $\SU(3)$-structure, such a structure is
  explicitly given in a standard basis in the appendix. We remark that
  most examples are constructed exploiting the stable form formalism
  and with computer support. In the following, we prove the
  non-existence of half-flat $\SU(3)$-structures on the remaining Lie
  algebras.

  In most of the cases, the obstructions of section \ref{sectionobst}
  can be applied directly.
  \begin{lemma}
    The Lie algebra $\g = \solv_2 \op \bR \op \solv_2 \op \bR $ and
    all Lie algebras $\g = \g_1 \op \g_2$ with $\g_1$ solvable and
    $\g_2$ one of the algebras $\solv_3$, $\solv_{3,\mu}\,,$ {\small
      $0 < |\mu| \le 1$}, $\solv'_{3,\mu}\,,$ {\small $\mu >0$}, do
    not admit a half-flat $\SU(3)$-structure.
  \end{lemma}
  \begin{proof}
    We want to apply the obstruction established in Proposition
    \ref{obst} and, given any of the Lie algebras $\g$ in question, we
    define a decomposition
    \[ V = \mathrm{span} \{ \e^1,\,\f^1\} \; , \quad W = \mathrm{span}
    \{ \e^2,\e^3,\f^2,\f^3\}\;, \] in a standard basis of $\g^*$.  By
    Lemma \ref{obst_particular_situation} it suffices to show that
    this is a coherent splitting such that the restrictions $d_{|
      \Lambda^3 W}$ and $d_{| \Lambda^4 W}$ are injective. In fact,
    the coherence can be verified directly by comparing the conditions
    of Lemma \ref{splittings}, (i), with the standard bases of the
    solvable three-dimensional Lie algebras.

    If $\g_2$ is one of the algebras $\solv_3$, $\solv_{3,\mu}$,
    {\small $0 < |\mu| \le 1$} or $\solv'_{3,\mu}$, {\small $\mu >0$},
    the standard bases satisfy
    \[d\f^2 \ne 0, \: \nexists\, c \in \bR: d\f^3 = c \, d\f^2, \:
    d\f^{23} \ne 0.\] Thus, considering again that the exterior
    algebras $\L^* \g_i^*$ of the summands are $d$-invariant, the
    image
    \[ d(\L^3 W) = \mathrm{span} \{ d(\e^{23}\f^2),\, d(\e^{23}\f^3),
    \, d(\e^2\f^{23}), \, d(\e^3\f^{23}) \}\] is four-dimensional and
    the image
    \[d(\L^4 W) = \mathrm{span} \{ d(\e^{23}\f^{23}) \}\] is
    one-dimensional. The same restrictions are injective for $\g =
    \solv_2 \op \bR \op \solv_2 \op \bR $, since in this case
    $d\e^{23} \ne 0$ and $d\f^{23} \ne 0$. This finishes the proof.
  \end{proof}
  The obstruction theory cannot be applied directly to the two
  remaining Lie algebras, although they admit coherent splittings and
  we have to deal with them separately.
  \begin{lemma}
    \label{lemmy}
    The Lie algebra $\g = \h_3 \op \solv_2 \op \bR$ does not admit a
    half-flat $\SU(3)$-structure. Furthermore, there is no
    decomposition $\g^* = V \op W$ with two-dimensional $V$ satisfying
    the obstruction condition \eqref{h03}.
  \end{lemma}
  \begin{proof}
    We start by proving the second assertion. Let $W \subset \g^*$ be
    an arbitrary four-dimensional subspace. It suffices to show that
    there is always a closed three-form with non-zero projection on
    $\L^3 W$. If the projection of $W$ on one of the summands $\g_i^*$
    is surjective, a generator of $\L^3 \g_i^*$ is closed and has
    non-zero projection on $\L^3 W$. For dimensional reasons, the only
    remaining possibility is that both projections have
    two-dimensional image in $W$. However, since all two-forms in
    $\L^2\h_3^*$ are closed and the kernel of $d$ is two-dimensional
    on $\solv_2 \op \bR$, there is necessarily a closed three-form in
    $\L^2 \h_3^* \^ (\solv_2 \op \bR)^*$ with non-zero projection on
    $\L^3 W$. Therefore, the obstruction condition \eqref{h03} is
    never satisfied.

    However, we can prove that there is no half-flat
    $\SU(3)$-structure by refining the idea of the obstruction
    condition as follows. Suppose that $(\rho,\o)$ is a half-flat
    $\SU(3)$-structure, i.e.\ $\rho \in Z^3$ and
    $\sigma=\frac{1}{2}\o^2 \in Z^4$ and let $\{ \e_1, \dots, \f_3\}$
    denote a standard basis of $\h_3 \op \solv_2 \op \bR$. We claim
    that
    \[ \f^1 \^ J_{\rho}^* \f^1 \^ \sigma = 0 \] which suffices to
    prove the non-existence since $\f_1$ would be isotropic by
    \eqref{gonL1}. First of all, an easy calculation reveals that
    \[ \f_1 \^ \sigma \in \mathrm{span} \{ \f^1\e^{12}\f^{23} , \,
    \f^1\e^{123}\f^{3} \} \] for an arbitrary closed four-form $\sigma$.
    Thus, it remains to show that $J_{\rho}^* \f^1$ has no component
    along $\e^3$ and $\f^2$ or equivalently that
    \[ J_{\rho}^* \f^1(v) \phi(\rho) \stackrel{~(\ref{JonL1})}{=} \f^1
    \^ (v \hook \rho) \^ \rho \] vanishes for $v \in \{ \e_3, \f_2 \}$.
    This assertion is straightforward to verify for an arbitrary
    closed three-form $\rho$.
  \end{proof}

  For the last remaining Lie algebra, we apply a different argument.
  \begin{lemma}
    \label{21obst}
    The Lie algebra $\g =\solv_2 \op \bR \op \bR^3$ does not admit a
    closed stable form $\rho$ with $\lambda(\rho) < 0$, in particular
    it does not admit a half-flat $\SU(3)$-structure. Furthermore,
    there is no decomposition $\g^* = V \op W$ with two-dimensional
    $V$ which satisfies the obstruction condition \eqref{h03}.
  \end{lemma}
  \begin{proof}
    Suppose that $\rho$ is a closed stable form inducing a complex or
    a para-complex structure $J_\rho$. Let $\{ \e_1, \e_2 \}$ be a
    basis of $\solv_2$ such that $d\e^2 = \e^{21}$. Since $\rho$ is
    closed, there are a one-form $\beta \in (\bR^4)^*$, a two-form
    $\gamma \in \L^2(\bR^4)^*$ and a three-form $\delta \in
    \L^3(\bR^4)^*$, such that
    \[\rho = \e^{12} \^ \b + \e^{1} \^ \gamma + \delta.\] Therefore,
    we have
    \[ K_{\rho}(\e_2) = \kappa ((\e_2 \hook \rho) \wedge \rho) =
    \kappa (-\e^{1} \^ \b \^ \delta) \] with $\b \^ \delta \in
    \L^4(\bR^4)^*$. However, this implies that $J(\e_2)$ is
    proportional to $\e_2$ by \eqref{J} which is only possible if
    $\lambda(\rho)>0$ and the first assertion is proven.

    In order to prove the second assertion, it suffices to show that
    for every four-dimensional subspace $W \subset \g^*$, there is a
    closed three-form with non-zero projection on $\L^3W$. This
    follows immediately from the observation that $\dim (\ker d) = 5$
    which implies that $\dim (\ker d \cap W) \geq 3$ for every
    four-dimensional subspace $W$.
  \end{proof}
  The lemma finishes the proof of the theorem as all possible direct
  sums according to the classification of three-dimensional Lie
  algebras have been considered.
\end{proof}

We remark that the lemmas \ref{lemmy} and \ref{21obst} give two
examples of solvable Lie algebras which show that the condition of
\cite[Theorem 5]{C}, which characterises six-dimensional nilpotent Lie
algebras admitting a half-flat $\SU(3)$-structure, cannot be
generalised without further restrictions to solvable Lie algebras.

\section{Half-flat $\SU(1,2)$-structures on direct sums}
\label{su21}
In this section, we describe some interesting observations concerning
the existence of half-flat $\SU(p,q)$-structures, $p+q =3$, with
indefinite metrics on direct sums $\g = \g_1 \op \g_2$ of
three-dimensional Lie algebras. It suffices to consider
$\SU(1,2)$-structures after possibly multiplying the metric by minus
one.

First of all, the obstruction condition of Proposition \ref{obst} does
not apply since isotropic subspaces are of course possible for metrics
of signature $(2,4)$. For instance, the Lie algebra $\solv_2 \op \bR
\op \solv_2 \op \bR$ does admit a half-flat $\SU(1,2)$-structure but
no half-flat $\SU(3)$-structure. Indeed, the structure defined in the
standard basis by
\begin{eqnarray*}
  \rho &=& - \e^{123} - \e^{12}\f^{3} - \e^{12}\f^{2}+2 \:\e^{13}\f^{3}+\:\e^{2}\f^{12} - \e^{3}\f^{13}+\:\f^{123}, \\
  \omega &=& \e^{13} - \e^{1}\f^{2}+\:\e^{1}\f^{3}+\:\e^{2}\f^{3} - \f^{12} ,\\
  g &=& - \, ( \e^2 )^2 - 2 \, ( \f^3 )^2 +2 \, \e^1 \ccdot \e^3 +2 \, \e^1 \ccdot \f^2 +2 \, \e^1 \ccdot \f^3 -2 \, \e^2 \ccdot \f^3 +2 \, \e^3 \ccdot \f^1 +2 \, \f^1 \ccdot \f^3,
\end{eqnarray*}
is a half-flat $\SU(1,2)$-structure with $V=\mathrm{span} \{ \e^1,\f^1 \}$
$J_\rho$-invariant and isotropic.

In fact, the obstruction established in Lemma \ref{21obst} is stronger
and also shows the non-existence of a half-flat $\SU(1,2)$-structure
on $\g =\solv_2 \op \bR \op \bR^3$. It can be generalised to the
following Lie algebras.
\begin{proposition}
  \label{21}
  Let $\g = \g_1 \op \g_2$ be a Lie algebra such that $\g_1$ is one of
  the algebras $\bR^3$, $\h_3$ or $\solv_2 \op \bR$ and $\g_2$ is one
  of the algebras $\solv_3$, $\solv_{3,\mu}\,,$ {\small $0 < |\mu| \le
    1$}, $\solv'_{3,\mu}\,,$ {\small $\mu >0$}.

  Every closed three-form $\rho$ on one of these Lie algebras $\g$
  satisfies $\lambda(\rho) \ge 0$. In particular, these Lie algebras
  do not admit a half-flat $\SU(p,q)$-structure for any signature
  $(p,q)$ with $p+q=3$.
\end{proposition}
\begin{proof}
  The proof is straightforward, but tedious without computer
  support. In a fixed basis, the condition $d\rho=0$ is linear in the
  coefficients of an arbitrary three-form $\rho$ and can be solved
  directly. When identifying $\L^6 V^*$ with $\bR$ with the help of a
  volume form $\nu$, one can calculate the quartic invariant
  $\lambda(\rho) \in \bR$, for instance in a standard basis. Carrying
  this out with Maple and factorising the resulting expression, we
  verified $\lambda(\rho) \ge 0$ for an arbitrary closed three-form on
  any of the Lie algebras in question.

  As a half-flat $\SU(p,q)$-structure is defined by a pair $(\rho,\o)$
  of stable forms which satisfy in particular $\la(\rho)<0$ and
  $d\rho=0$, such a structure cannot exist and the lemma is proven.
\end{proof}

We add the remark, that a result analogous to Lemma $\ref{u3ortho}$
for a pseudo-Hermitian structure of indefinite signature would involve
a considerably more complicated normal form for $\omega$. Therefore, a
generalisation of the proof of Theorem \ref{theo_orth} to indefinite
metrics seems to be difficult.

\section{Half-flat $\SL(3,\bR)$-structures on direct sums}
\label{sl3}
Finally, we turn to the para-complex case of $\SL(3,\bR)$-structures.
As explained in section \ref{prelim}, a half-flat
$\SL(3,\bR)$-structure is defined by a pair $(\rho,\o)$ of stable
forms such that $J_\rho$ is an almost para-complex structure and
\begin{equation*}
  \o \wedge \rho = 0, \quad d\o^2 = 0, \quad d\rho=0 .  
\end{equation*}
As the induced metric is always neutral and $\la(\rho)>0$, neither
Proposition \ref{obst} nor Lemma \ref{21} obstruct the existence of
such a structure. For instance, the Lie algebra $\solv_2 \op \bR \op
\solv_3$ does not admit a half-flat $\SU(p,q)$-structure for any
signature $(p,q)$ with $p+q=3$ , but
\begin{eqnarray*}
  \rho &=& -2 \:\e^{12}\f^{3}-2 \:\e^{2}\f^{31}+\:\e^{3}\f^{12} - \e^{3}\f^{31}+\:\f^{123}, \\
  \o &=& \e^{13} - \e^{23} +\:\e^{1}\f^{3} +\:\e^{2}\f^{2}  - \e^{3}\f^{1} +2 \:\f^{13},  \\
  g &=& -2 \, ( \, \e^1 \ccdot \e^3  - \, \e^2 \ccdot \e^3 + \, \e^1 \ccdot \f^3 + \, \e^2 \ccdot \f^2 + \, \e^3 \ccdot \f^1 \, ),
\end{eqnarray*}
is an example of a half-flat $\SL(3,\bR)$-structure.

When trying to generalise Theorem \ref{theo_orth} to the para-complex
situation, we find an astonishingly similar result if we additionally
require the metric to be definite when restricted to one of the
summands. We omit the proofs which are very similar to the original
ones due to the analogies explained in section \ref{prelim}.

\begin{lemma}
  Let $(V_1,g_1)$ and $(V_2,g_2)$ be Euclidean vector spaces and let
  $(g,J,\omega)$ be a para-Hermitian structure on the orthogonal
  product $(V_1 \oplus V_2, g = - g_1 + g_2)$. There are orthonormal
  bases $\{ e_1, e_2, e_3 \}$ of $V_1$ and $\{ f_1, f_2, f_3 \}$ of
  $V_2$ which can be joined to a pseudo-orthonormal basis of $V_1
  \oplus V_2$ such that
  \begin{equation}
    \omega = a \, \e^{12} + \sqrt{1 + a^2} \, \e^{1} \f^1 + \sqrt{1 + a^2} \, \e^2 \f^2 + \e^3 \f^3 + a \, \f^{12}
  \end{equation}
  for a real number $a$.
\end{lemma}
In analogy to the Hermitian case, we call the para-Hermitian structure
of type I if it admits a basis with $a=0$ and of type II if it admits
a basis with $a\ne 0$.

\begin{theorem}
  Let $\g=\g_1 \op \g_2$ be a direct sum of three-dimensional Lie
  algebras $\g_1$ and $\g_2$. 

  The Lie algebra $\g$ admits a half-flat $\SL(3,\bR)$-structure
  such that $\g_1$ and $\g_2$ are mutually orthogonal, such that the
  restriction of the metric to both summands is definite and such that
  the underlying para-Hermitian structure is of type I if and only if
  \begin{align*}
    \mbox{(i)} & \quad \mbox{$\g_1=\g_2$ and both are unimodular or } \\
    \mbox{(ii)} & \quad \mbox{$\g_1$ is non-abelian unimodular and $\g_2$ abelian or vice versa.} 
  \end{align*}
  Moreover, the Lie algebra $\g$ admits a half-flat
  $\SL(3,\bR)$-structure such that $\g_1$ and $\g_2$ are mutually
  orthogonal, such that the restriction of the metric to both summands
  is definite and such that the underlying para-Hermitian structure is
  of type II if and only if the pair $(\g_1,\g_2)$ or $(\g_2,\g_1)$ is
  contained in the following list:
  \begin{align*}
    &&& (\mathfrak e(1,1),\mathfrak e(1,1)), \\
    &&& (\mathfrak e(2),\bR \oplus \solv_{2}), \\
    &&& (\su(2),\solv_{3,\mu})  && \mbox{for $0 < \mu \le 1$,} && \\
    &&& (\sl(2,\bR),\solv_{3,\mu}) && \mbox{for $-1 < \mu < 0$.} &&
  \end{align*}
\end{theorem}
If we require, instead of orthogonality, that the
$\SL(3,\bR)$-structure is adapted to the direct sum $\g_1 \oplus
{\g_2}$ in the sense that the summands $\g_1$ and $\g_2$ are the
eigenspaces of $J$, we find the following interesting relation to
unimodularity.
\begin{proposition}
  \label{para2}
  A direct sum $\g_1 \op \g_2$ of three-dimensional Lie algebras
  $\g_1$ and $\g_2$ admits a half-flat $\SL(3,\bR)$-structure
  $(g,J,\omega,\Psi)$ such that $\g_1$ and ${\g_2}$ are the $\pm
  1$-eigenspaces of $J$ if and only if both $\g_1$ and $\g_2$ are
  unimodular.
\end{proposition}
\begin{proof}
  Let $(g,J,\omega,\Psi)$ be an $\SL(3,\bR)$-structure on $\g=\g_1 \op
  \g_2$ such that $\g_1$ is the $+1$-eigenspace of $J$ and $\g_2$ is
  the $-1$-eigenspaces of $J$. Since $\psi^+ = \Re \Psi$ is a stable
  form inducing the para-complex structure $J$, we can choose bases
  $\{ \e^i \}$ of $\g_1^*$ and $\{ \f^i \}$ of $\g_2^*$ such that
  $\psi^+ = \e^{123} + \f^{123}$ is in the normal form
  \ref{paranormal}. Thus, the real part $\psi^+$ is closed as we are
  dealing with a direct sum of Lie algebras. Due to the simple form of
  $\psi^+$, it is easy to verify that the relation $\omega \wedge
  \psi^+ = 0$ holds for an arbitrary non-degenerate $\omega$ if and
  only if $\omega$ has only terms in $\g_1^* \otimes \g_2^*$. Now
  we are in the situation of Lemma \ref{omegasquare} and conclude that
  the only remaining equation $d \omega^2 =0$ is satisfied if and only
  if both $\g_1$ and $\g_2$ are unimodular.
\end{proof}

\clearpage

\section{Appendix}
The following tables contain all direct sums $\g =\g_1 \op \g_2 $ of
three-dimensional Lie algebras which admit a half-flat
$\SU(3)$-structure. In each case, an explicit example $(\o,\rho)$ of a
normalised half-flat $\SU(3)$-structure including the induced
Riemannian metric $g$ is given where $\{ \e^i \}$ is a standard basis
(as defined in tables \ref{unilist} and \ref{nonunilist}) of
$\g_1^*$ and $\{ \f^i \}$ is a standard basis of $\g_2^*$.
\vspace{0.05cm} {\small
  \begin{table}[ht]
    \caption{Unimodular direct sums of three-dimensional Lie algebras}
    \label{uni}
    \renewcommand{\arraystretch}{1.5}
    \setlength{\tabcolsep}{0.3cm}
    \begin{tabular}{|c|l|}
      \hline { Lie algebra} & { Half-flat $\SU(3)$-structure with $\o = \e^1\f^1 + \e^2\f^2 + \e^3\f^3$ }  \\ \hline 
      $\h \op \h$, & 
      $\rho = \frac{1}{2}\sqrt 2 \,\{\, \e^{123} - \e^{1}\f^{23} - \e^{2}\f^{31} - \e^{3}\f^{12} + \e^{12}\f^{3} + \e^{31}\f^{2} + \e^{23}\f^{1} - \f^{123} \,\}$ \\ 
      $\h$ unimodular
      & $g= ( \e^1 )^2 + ( \e^2 )^2 + ( \e^3 )^2 + ( \f^1 )^2 + ( \f^2 )^2 + ( \f^3 )^2$ \\\hline 
      $\h \op \bR^3$, &
      $\rho = \e^{12}\f^{3} + \e^{31}\f^{2} + \e^{23}\f^{1} - \f^{123}$ \\
      $\h$ unimodular
      & $g= ( \e^1 )^2 + ( \e^2 )^2 + ( \e^3 )^2 + ( \f^1 )^2 + ( \f^2 )^2 + ( \f^3 )^2$ \\ \hline
      $\su(2) \op \sl(2,\bR)$, & 
      $\rho = 2^\frac{1}{4} \,\{\, \frac{1}{2}\,\e^{123} + \e^{23}\f^{1} + \e^{31}\f^{2} + \e^{12}\f^{3} - \e^{1}\f^{23} - \e^{2}\f^{31} + \e^{3}\f^{12}-2\,\f^{123} \,\}$ \\
      & $g= \sqrt 2 \,\{\, \frac{3}{2} \, ( \e^1 )^2 + \frac{3}{2} \, ( \e^2 )^2 +\frac{1}{2} \, ( \e^3 )^2 + ( \f^1 )^2 + ( \f^2 )^2 +3 \, ( \f^3 )^2$ \\ 
      & $\quad + \, 2 \, \e^1 \ccdot \f^1 +2 \, \e^2 \ccdot \f^2 -2 \, \e^3 \ccdot \f^3 \,\}$ \\ \hline

      $\su(2) \op \efr(2)$ & 
      $\rho =  - \e^{23}\f^{1} - \e^{31}\f^{2} - \e^{12}\f^{3} + \e^{2}\f^{31} + \e^{3}\f^{12} + \f^{123}$ \\ 
      & $g = \, ( \e^1 )^2 + ( \e^2 )^2 + ( \e^3 )^2 +2 \, ( \f^1 )^2 + ( \f^2 )^2 + ( \f^3 )^2 -2 \, \e^1 \ccdot \f^1$  \\ \hline

      $ \sl(2,\bR) \op \efr(2)$ &
      $\rho = -2 \,\e^{23}\f^{1} - \e^{31}\f^{2} - \e^{12}\f^{3} + \e^{2}\f^{31} - \e^{3}\f^{12} + \f^{123}$ \\ &
      $g = \, ( \e^1 )^2 +2 \, ( \e^2 )^2 +2 \, ( \e^3 )^2 + ( \f^1 )^2 + ( \f^2 )^2 + ( \f^3 )^2 +2 \, \e^2 \ccdot \f^2 -2 \, \e^3 \ccdot \f^3 $ \\ \hline

      $ \su(2) \op \efr(1,1)$, &
      $\rho = -2 \,\e^{23}\f^{1} - \e^{31}\f^{2} - \e^{12}\f^{3} + \e^{2}\f^{31} - \e^{3}\f^{12} + \f^{123}$ \\
      $\efr(2) \op \efr(1,1)$ & $g = \, ( \e^1 )^2 +2 \, ( \e^2 )^2 +2 \, ( \e^3 )^2 + ( \f^1 )^2 + ( \f^2 )^2 + ( \f^3 )^2 +2 \, \e^2 \ccdot \f^2 -2 \, \e^3 \ccdot \f^3$  \\ \hline

      $ \sl(2,\bR) \op \efr(1,1)$ &
      $\rho =  - \e^{23}\f^{1} - \e^{31}\f^{2} - \e^{12}\f^{3} + \e^{2}\f^{31} + \e^{3}\f^{12} + \f^{123}$ \\
      & $g = \, ( \e^1 )^2 + ( \e^2 )^2 + ( \e^3 )^2 +2 \, ( \f^1 )^2 + ( \f^2 )^2 + ( \f^3 )^2 -2 \, \e^1 \ccdot \f^1$  \\ \hline

      $\su(2) \op \h_3$, & 
      $ \rho =  - \e^{23}\f^{1}-\frac{5}{4} \,\e^{31}\f^{2} - \e^{12}\f^{3} + \e^{3}\f^{12} + \f^{123}$  \\
      $\efr(2) \op \h_3$ & 
      $g= \frac{5}{4} \, ( \e^1 )^2 + ( \e^2 )^2 +\frac{5}{4} \, ( \e^3 )^2 + ( \f^1 )^2 +\frac{5}{4} \, ( \f^2 )^2 + ( \f^3 )^2$ \\
      & $\quad - \e^1 \ccdot \f^1 - \e^2 \ccdot \f^2 + \e^3 \ccdot \f^3$ \\ \hline

      $ \sl(2,\bR) \op \h_3$,  &
      $\rho =  - \e^{23}\f^{1}-\frac{5}{4} \,\e^{31}\f^{2} - \e^{12}\f^{3} - \e^{3}\f^{12} + \f^{123}$ \\
      $\efr(1,1) \op \h_3 $ & 
      $ g= \frac{5}{4} \, ( \e^1 )^2 + ( \e^2 )^2 +\frac{5}{4} \, ( \e^3 )^2 + ( \f^1 )^2 +\frac{5}{4} \, ( \f^2 )^2 + ( \f^3 )^2$ \\
      & $\quad +\, \e^1 \ccdot \f^1 + \e^2 \ccdot \f^2 - \e^3 \ccdot \f^3 $ \\ \hline
    \end{tabular}
  \end{table}

  \begin{table}[ht]
    \caption{Solvable, non-unimodular direct sums admitting a half-flat $\SU(3)$-structure}
    \label{nonuni}
    \renewcommand{\arraystretch}{1.5}
    \setlength{\tabcolsep}{0.3cm}
    \begin{tabular}{|c|l|}
      \hline { Lie algebra} & { Half-flat $\SU(3)$-structure } \\ \hline
      $\efr(2) \op \solv_2 \op \bR$ &
      $\o =  \e^{12} + \e^{3}\f^{1} - \f^{23} $ \\
      & $\rho = \e^{23}\f^{3} + \e^{2}\f^{21} + \e^{13}\f^{2} - \e^{1}\f^{31} $ \\
      & $g = ( \e^1 )^2 + ( \e^2 )^2 + ( \e^3 )^2 + ( \f^1 )^2 + ( \f^2 )^2 + ( \f^3 )^2 $ \\ \hline
      $\efr(1,1) \op \solv_2 \op \bR$ &
      $\o =  - \e^{1}\f^{3} - \e^{3}\f^{2} + \e^{2}\f^{1}  - \f^{23}$ \\
      & $\rho = \e^{23}\f^{3}-2 \,\e^{31}\f^{1} + \e^{12}\f^{2}-3 \,\e^{1}\f^{31} - \e^{3}\f^{12}+2 \,\f^{123} $ \\
      & $g = 2 \, ( \e^1 )^2 + ( \e^2 )^2 +2 \, ( \e^3 )^2 + ( \f^1 )^2 + ( \f^2 )^2 +5 \, ( \f^3 )^2 - 2 \, \e^1 \ccdot \f^2 -6 \, \e^3 \ccdot \f^3$ \\ \hline
    \end{tabular}
  \end{table}

  \begin{table}[ht]
    \caption{Direct sums which are neither solvable nor unimodular}
    \label{nonsolv}
    \renewcommand{\arraystretch}{1.5}
    \setlength{\tabcolsep}{0.3cm}
    \begin{tabular}{|c|l|}
      \hline
      { Lie algebra} & { Half-flat $\SU(3)$-structure } \\ \hline
      $\su(2) \op \solv_2 \op \bR$,
      & $\o =  \e^{1}\f^{1} - \f^{23} + \e^{2}\f^{2} + \e^{3}\f^{3}$ \\
      $\sl(2,\bR) \op \solv_2 \op \bR$
      & $\rho = \e^{23}\f^{1} + \e^{31}\f^{2} + \e^{12}\f^{3} + \e^{2}\f^{12} - \f^{123}$ \\
      & $g = \, ( \e^1 )^2 + ( \e^2 )^2 + ( \e^3 )^2 + ( \f^1 )^2 +2 \, ( \f^2 )^2 + ( \f^3 )^2 -2 \, \e^3 \ccdot \f^2$ \\ \hline
      $\su(2) \op \solv_{3}$ & 
      $\omega = \f^{23} + \e^{23}+2 \,\e^{1}\f^{1}$ \\ &
      $\rho = \frac{2}{3}\, 3^\frac{3}{4} \, \{ \,\e^{31}\f^{2} - \e^{12}\f^{3} - \e^{2}\f^{31} + \e^{3}\f^{31} + \e^{2}\f^{12} \, \}  $ \\ &
      $g = \frac{2}{3}\sqrt 3 \, \{ \, 2 \, ( \e^1 )^2  +  ( \e^2 )^2 + ( \e^3 )^2 + 2 \, ( \f^1 )^2 + ( \f^2 )^2 + ( \f^3 )^2 $\\
      &$\quad + \,2 \, \e^1 \ccdot \f^1 - \e^2 \ccdot \e^3 + \f^2 \ccdot \f^3 \}$ \\ \hline
      $\sl(2) \op \solv_{3}$ & 
      $\omega = \e^{1}\f^{1}-2 \,\f^{23} + \e^{3}\f^{3} + \e^{2}\f^{2}$ \\ 
      & $\rho = \frac{1}{3} \,\e^{23}\f^{1}+3 \,\e^{31}\f^{2} + \e^{31}\f^{3} + \e^{12}\f^{2}+\frac{4}{3} \,\e^{12}\f^{3}-4 \,\e^{2}\f^{31}+\frac{7}{3} \,\e^{3}\f^{31}+3 \,\e^{2}\f^{12}$ \\ &
      $\quad  - \e^{3}\f^{12}-26 \,\f^{123} $ \\
      & $g= 3 \, ( \e^1 )^2 +\frac{4}{9} \, ( \e^2 )^2 + ( \e^3 )^2 +\frac{17}{3} \, ( \f^1 )^2 +94 \, ( \f^2 )^2 +\frac{328}{9} \, ( \f^3 )^2 $ \\ 
      & $\quad - 8 \, \e^1 \ccdot \f^1 -\frac{2}{3} \, \e^2 \ccdot \e^3 +\frac{34}{3} \, \e^2 \ccdot \f^2 +\frac{16}{9} \, \e^2 \ccdot \f^3 -16 \, \e^3 \ccdot \f^2 -\frac{34}{3} \, \e^3 \ccdot \f^3 +\frac{224}{3} \, \f^2 \ccdot \f^3$ \\ \hline

      $\su(2) \op \solv_{3,\mu}$ &
      $\omega = \frac{1}{\mu+1} \,\e^{12} + \e^{3}\f^{1} - \f^{32} $ \\ 
      {\small $( 0 < \mu \le 1)$} &
      $\rho = \mu^{-\frac{1}{4}} (\mu+1)^{-\frac{1}{2}}\, \{ \,\e^{13}\f^{2} - \e^{23}\f^{3}-\mu \,\e^{1}\f^{13} - \e^{2}\f^{12} \, \} $ \\ &
      $g = \mu^{-\frac{1}{2}} \, \{ \, \frac{\mu}{\mu+1} \, ( \e^1 )^2 + \frac{1}{\mu+1} \, ( \e^2 )^2 +  \, ( \e^3 )^2 + \mu \, ( \f^1 )^2 + \, ( \f^2 )^2 + \mu \, ( \f^3 )^2 \, \} $ \\ \hline

      $\sl(2) \op \solv_{3,\mu}$ &
      $\omega = \frac{1}{\mu+1} \,\e^{23} + \e^{1}\f^{1} + \f^{32} $ \\ 
      {\small $( -1 < \mu < 0)$ } &
      $\rho = (-\mu)^{-\frac{1}{4}} (\mu+1)^{-\frac{1}{2}} \, \{  \,\e^{12}\f^{3}  - \e^{13}\f^{2}  + \e^{2}\f^{12} - \mu \,\e^{3}\f^{13} \, \} $ \\ &
      $g =(-\mu)^{-\frac{1}{2}} \, \{ \, ( \e^1 )^2 +\frac{1}{\mu+1} \, ( \e^2 )^2 - \frac{\mu}{\mu+1} \, ( \e^3 )^2 -\mu \, ( \f^1 )^2 + ( \f^2 )^2 -\mu \, ( \f^3 )^2 \, \}$ \\ \hline

      $\su(2) \op \solv_{3,\mu}$ &
      $\omega =  \f^{23} + \e^{3}\f^{1}-\frac{\mu (2 \mu+3)}{2(\mu+1)^2} \,\e^{23} - \e^{1}\f^{1} + \e^{1}\f^{3}+\frac{\mu (2 \mu+3)}{2(\mu+1)^2} \,\e^{12}- \frac{2\mu^2 +\mu -2}{2(\mu+1)^2} \,\e^{2}\f^{2} + \e^{3}\f^{3} $ \\ 
      {\small $( -1 < \mu < 0)$ } & 
      $\rho = -\frac{2\mu^2+3\mu+2}{2(\mu+1)^2} \,\e^{23}\f^{1}- \frac{1}{\mu} \,\e^{23}\f^{3}-2 \,\e^{13}\f^{2}+\frac{2\mu^2+3\mu+2}{2(\mu+1)^2} \,\e^{12}\f^{1} $ \\ &
      $\quad - \, \frac{1}{\mu} \,\e^{12}\f^{3} - \e^{1}\f^{13} - \e^{3}\f^{13}+2 \,\e^{2}\f^{12}+2 \,\f^{123} $ \\ &
      $g= -\frac{\mu^2+\mu+1}{\mu(\mu+1)} \, ( \e^1 )^2 -\frac{4\mu^4 +20\mu^3 + 29\mu^2+16\mu+4}{4 \mu (\mu+1)^3} \, ( \e^2 )^2 -\frac{\mu^2+\mu+1}{\mu(\mu+1)} \, ( \e^3 )^2 $\\
      &$\quad - \frac{\mu}{\mu+1} \, ( \f^1 )^2  +  \frac{4+3\mu}{\mu+1} \, ( \f^2 )^2 
      - \frac{\mu+1}{\mu} \, ( \f^3 )^2$\\
      &$\quad + \frac{ 2(\mu^2+1+3\mu)}{\mu(\mu+1)} \, \e^1 \ccdot \e^3 +\frac{2(\mu+2)}{\mu+1} \, \e^1 \ccdot \f^2 
      -\frac{2\mu^2 +5\mu +2}{\mu(\mu+1)} \, \e^2 \ccdot \f^3 +\frac{2(\mu+2)}{\mu+1} \, \e^3 \ccdot \f^2 $ \\ \hline

      $\sl(2) \op \solv_{3,\mu}$  &
      $\omega = \frac{2 (2 \mu+1)^{\frac{1}{2}}}{(\mu+1)^2} \,\e^{1}\f^{3}  + \e^{2}\f^{1} + \f^{23}+ \frac{\mu}{\mu+1} \,\e^{13} + \e^{1}\f^{2} + \e^{3}\f^{3} $ \\ 
      {\small $( 0 < \mu \le 1)$ } &
      $\rho =2 \frac{2 (2 \mu+1)^{\frac{1}{2}}}{(\mu+1)^2}\,\e^{123}
      + \e^{23}\f^{2} - \e^{13}\f^{1}+\frac{1}{\mu}\,\e^{12}\f^{3} - \e^{3}\f^{13} + \e^{1}\f^{12}+ \frac{\mu+1}{\mu} \,\f^{123} $ \\ &
      $g = \frac{\mu^3+11\mu^2+7\mu+1}{\mu(\mu+1)^3} \, ( \e^1 )^2 +\frac{\mu+1}{\mu} \, ( \e^2 )^2 + (2\mu+1) \, ( \e^3 )^2 + \frac{\mu+1}{\mu} \, ( \f^1 )^2 +\frac{\mu+1}{\mu^2} \, ( \f^3 )^2 $ \\ &
      $\quad +\frac{1+3\mu+2\mu^2}{\mu} \, ( \f^2 )^2 + \frac{6(2\mu+1)^{\frac{1}{2}}}{\mu+1} \, \e^1 \ccdot \e^3 +\frac{2(2\mu+1)^{\frac{1}{2}}(3\mu+1)}{\mu(\mu+1)} \, \e^1 \ccdot \f^2 + \frac{4(2\mu+1)}{\mu(\mu+1)^2} \, \e^1 \ccdot \f^3 $ \\ &
      $\quad +\frac{2(2\mu+1)^{\frac{1}{2}}}{\mu} \, \e^2 \ccdot \f^1 + (4+4\mu) \, \e^3 \ccdot \f^2 +\frac{2(2\mu+1)^{\frac{1}{2}}}{\mu} \, \e^3 \ccdot \f^3 +\frac{2(2\mu+1)^{\frac{1}{2}}}{\mu} \, \f^2 \ccdot \f^3 $ \\ \hline

      $\su(2) \op \solv'_{3,\mu}$ &
      $\o =  \e^{2}\f^{2}-2 \mu \,\f^{23} + \e^{3}\f^{3} + \e^{1}\f^{1} $ \\ 
      {\small ($\mu > 0$)} &
      $\rho = \,\e^{23}\f^{1} + \e^{31}\f^{2} + \e^{12}\f^{3} + \e^{2}\f^{31}-\mu \,\e^{3}\f^{31}+\mu \,\e^{2}\f^{12} + \e^{3}\f^{12}+(\mu^2-1) \,\f^{123} $ \\ &
      $g =  \, ( \e^1 )^2 + ( \e^2 )^2 + ( \e^3 )^2 +2 \, ( \f^1 )^2 +(\mu^2+1) \, ( \f^2 )^2 +(\mu^2+1) \, ( \f^3 )^2 $ \\ & 
      $ \quad + \, 2 \, \e^1 \ccdot \f^1 +2\mu \, \e^2 \ccdot \f^3 -2\mu \, \e^3 \ccdot \f^2$ \\ \hline

      $\sl(2,\bR) \op \solv'_{3,\mu}$ &
      $\omega =  \e^{2}\f^{2}-2 \mu \,\f^{23} + \e^{3}\f^{3} + \e^{1}\f^{1} $ \\ 
      {\small ($\mu > 0$)} &
      $\rho= \frac{1}{2} \,\e^{23}\f^{1}+2 \,\e^{31}\f^{2} + \e^{12}\f^{3}+2 \,\e^{2}\f^{31}+\mu \,\e^{3}\f^{31}$\\
      &$\quad +\,2 \mu \,\e^{2}\f^{12} - \e^{3}\f^{12}- (4 \mu^2+ \frac{29}{4}) \,\f^{123}  $ \\ &
      $g  = 2 \, ( \e^1 )^2 +\frac{1}{2} \, ( \e^2 )^2 + ( \e^3 )^2 + \frac{13}{8} \, ( \f^1 )^2 +(16 \mu^2 + \frac{29}{2}) \, ( \f^2 )^2 +(2 \mu^2 + \frac{29}{4}) \, ( \f^3 )^2 $ \\ & 
      $ \quad + \,  3 \, \e^1 \ccdot \f^1 -5 \, \e^2 \ccdot \f^2 -2\mu \, \e^2 \ccdot \f^3 -8\mu \, \e^3 \ccdot \f^2 +5 \, \e^3 \ccdot \f^3 -10\mu \, \f^2 \ccdot \f^3$ \\\hline

    \end{tabular}
  \end{table}
}

\end{document}